\documentclass[twoside]{article}

\usepackage[accepted]{aistats2021}
\usepackage[round]{natbib}

\usepackage[colorlinks,
	linkcolor=blue,
	citecolor=blue,
	urlcolor=magenta,
	linktocpage,
	plainpages=false]{hyperref}

\usepackage{amsmath, amssymb, amsthm, graphicx,url,multicol}
\usepackage[dvipsnames]{xcolor}
\usepackage{algorithm}
\usepackage{algorithmic}

\usepackage{enumitem}

\usepackage {tikz}
\usepackage{pgfkeys}
\usetikzlibrary{shapes,arrows}
\usepackage{pgfplots}
\pgfplotsset{compat=1.13}
\pgfplotsset{plotOptions/.style={%
		label style={font=\scriptsize},
		legend style={font=\scriptsize},
		tick label style={font=\scriptsize},
		solid,
		very thick
	}}
\definecolor{colorP1}{RGB}{55,126,184}  
\definecolor{colorP2}{RGB}{228,26,28}  
\definecolor{colorP3}{RGB}{152,78,163} 
\definecolor{colorP4}{RGB}{77,175,74}  
\definecolor{colorP5}{RGB}{250, 150, 10} 
\definecolor{colorP6}{RGB}{139,69,19}

\renewcommand{\algorithmicrequire}{\textbf{Input:}} 
\renewcommand{\algorithmicensure}{\textbf{Output:}}

\newcommand{\N}{\mathbb{N}}
\newcommand{\lnorm}{\left\|}
\newcommand{\rnorm}{\right\|}

\begin{document}

\definecolor{ddarkbrown}{rgb}{0.5,0.2,0.05} \definecolor{bbluegray}{rgb}{0.05,0,0.5}

\newtheorem{theorem}{Theorem}[section]
\newtheorem{proposition}[theorem]{Proposition}
\newtheorem{definition}[theorem]{Definition}
\newtheorem{lemma}[theorem]{Lemma}
\newtheorem{corollary}[theorem]{Corollary}
\renewenvironment{proof}{\textbf{Proof.}}{\QED\bigskip}
\newtheorem{assumption}[theorem]{Assumption}
\newtheorem{remark}[theorem]{Remark}

\def\xx{{\boldsymbol x}}
\def\yy{{\boldsymbol y}}
\def\vv{{\boldsymbol v}}

\newcommand{\BEAS}{\begin{eqnarray*}}
\newcommand{\EEAS}{\end{eqnarray*}}
\newcommand{\BEA}{\begin{eqnarray}}
\newcommand{\EEA}{\end{eqnarray}}
\newcommand{\BEQ}{\begin{equation}}
\newcommand{\EEQ}{\end{equation}}
\newcommand{\BEQS}{\begin{equation*}}
\newcommand{\EEQS}{\end{equation*}}
\newcommand{\BIT}{\begin{itemize}}
\newcommand{\EIT}{\end{itemize}}
\newcommand{\BNUM}{\begin{enumerate}}
\newcommand{\ENUM}{\end{enumerate}}

\newcommand{\BA}{\begin{array}}
\newcommand{\EA}{\end{array}}

\newcommand{\refp}[1]{(\ref{#1})}

\newcommand{\cf}{{\it cf.}}
\newcommand{\eg}{{\it e.g.}}
\newcommand{\ie}{{\it i.e.}}
\newcommand{\etc}{{\it etc.}}
\newcommand{\ones}{\mathbf 1}

\newcommand{\reals}{{\mathbb R}}
\newcommand{\sreals}{\scriptsize{\mbox{\bf R}}}
\newcommand{\integers}{{\mbox{\bf Z}}}
\newcommand{\eqbydef}{\mathrel{\stackrel{\Delta}{=}}}
\newcommand{\complex}{{\mbox{\bf C}}}
\newcommand{\symm}{{\mbox{\bf S}}}  

\newcommand{\Span}{\mbox{\textrm{span}}}
\newcommand{\Range}{\mbox{\textrm{range}}}
\newcommand{\nullspace}{{\mathcal N}}
\newcommand{\range}{{\mathcal R}}
\newcommand{\diam}{\mathop{\bf diam}}
\newcommand{\sphere}{{\mathbb S}}
\newcommand{\Nullspace}{\mbox{\textrm{nullspace}}}
\newcommand{\Rank}{\mathop{\bf Rank}}
\newcommand{\NumRank}{\mathop{\bf NumRank}}
\newcommand{\NumCard}{\mathop{\bf NumCard}}
\newcommand{\Card}{\mathop{\bf Card}}
\newcommand{\Tr}{\mathop{\bf Tr}}
\newcommand{\diag}{\mathop{\bf diag}}
\newcommand{\lambdamax}{{\lambda_{\rm max}}}
\newcommand{\lambdamin}{\lambda_{\rm min}}
\newcommand{\idm}{\mathbf{I}}

\newcommand{\Expect}{\textstyle\mathop{\bf E}}
\newcommand{\Median}{\textstyle\mathop{\bf M}}
\newcommand{\Prob}{\mathop{\bf Prob}}
\newcommand{\erf}{\mathop{\bf erf}}

\newcommand{\Co}{{\mathop {\bf Co}}}
\newcommand{\co}{{\mathop {\bf Co}}}
\newcommand{\Po}{{\mathop {\bf Po}}}
\newcommand{\Ext}{{\mathop {\bf Ext}}}
\newcommand{\Var}{\mathop{\bf var{}}}
\newcommand{\dist}{\mathop{\bf dist{}}}
\newcommand{\Ltwo}{{\bf L}_2}
\newcommand{\QED}{~~\rule[-1pt]{6pt}{6pt}}\def\qed{\QED}
\newcommand{\approxleq}{\mathrel{\smash{\makebox[0pt][l]{\raisebox{-3.4pt}{\small$\sim$}}}{\raisebox{1.1pt}{$<$}}}}
\newcommand{\argmin}{\mathop{\rm argmin}}
\newcommand{\epi}{\mathop{\bf epi}}
\newcommand{\var}{\mathop{\bf var}}

\newcommand{\vol}{\mathop{\bf vol}}
\newcommand{\Vol}{\mathop{\bf vol}}

\newcommand{\dom}{\mathop{\bf dom}}
\newcommand{\aff}{\mathop{\bf aff}}
\newcommand{\cl}{\mathop{\bf cl}}
\newcommand{\Angle}{\mathop{\bf angle}}
\newcommand{\intr}{\mathop{\bf int}}
\newcommand{\relint}{\mathop{\bf rel int}}
\newcommand{\bd}{\mathop{\bf bd}}
\newcommand{\vect}{\mathop{\bf vec}}
\newcommand{\dsp}{\displaystyle}
\newcommand{\foequal}{\simeq}
\newcommand{\VOL}{{\mbox{\bf vol}}}
\newcommand{\argmax}{\mathop{\rm argmax}}
\newcommand{\xopt}{x^{\rm opt}}

\newcommand{\Xb}{{\mbox{\bf X}}}
\newcommand{\xst}{x^\star}
\newcommand{\varphist}{\varphi^\star}
\newcommand{\lambdast}{\lambda^\star}
\newcommand{\Zst}{Z^\star}
\newcommand{\fstar}{f^\star}
\newcommand{\xstar}{x^\star}
\newcommand{\xc}{x^\star}
\newcommand{\lambdac}{\lambda^\star}
\newcommand{\lambdaopt}{\lambda^{\rm opt}}

\newcommand{\geqK}{\mathrel{\succeq_K}}
\newcommand{\gK}{\mathrel{\succ_K}}
\newcommand{\leqK}{\mathrel{\preceq_K}}
\newcommand{\lK}{\mathrel{\prec_K}}
\newcommand{\geqKst}{\mathrel{\succeq_{K^*}}}
\newcommand{\gKst}{\mathrel{\succ_{K^*}}}
\newcommand{\leqKst}{\mathrel{\preceq_{K^*}}}
\newcommand{\lKst}{\mathrel{\prec_{K^*}}}
\newcommand{\geqL}{\mathrel{\succeq_L}}
\newcommand{\gL}{\mathrel{\succ_L}}
\newcommand{\leqL}{\mathrel{\preceq_L}}
\newcommand{\lL}{\mathrel{\prec_L}}
\newcommand{\geqLst}{\mathrel{\succeq_{L^*}}}
\newcommand{\gLst}{\mathrel{\succ_{L^*}}}
\newcommand{\leqLst}{\mathrel{\preceq_{L^*}}}
\newcommand{\lLst}{\mathrel{\prec_{L^*}}}

\newcommand{\realsp}{\mathbf{R}_+^n}
\newcommand{\intrealsp}{\int_{\mathbf{R}_+^n}}

%

%

\twocolumn[

\aistatstitle{Convergence of Constrained Anderson Acceleration}

\aistatsauthor{ Mathieu Barr\'e \And Adrien Taylor \And Alexandre d'Aspremont }

\aistatsaddress{ INRIA - D.I ENS, \\  ENS, CNRS, PSL University,\\ Paris, France \And INRIA - D.I ENS, \\  ENS, CNRS, PSL University,\\ Paris, France \And CNRS - D.I ENS, \\  ENS, CNRS, PSL University,\\ Paris, France } ]
\begin{abstract}
  We prove non asymptotic linear convergence rates for the constrained Anderson acceleration extrapolation scheme. These guarantees come from new upper bounds on the constrained Chebyshev problem, which consists in minimizing the maximum absolute value of a polynomial on a bounded real interval with $l_1$ constraints on its coefficients vector. Constrained Anderson Acceleration has a numerical cost comparable to that of the original scheme.
\end{abstract}

\section{Introduction}
Obtaining faster convergence rates is a central concern in numerical analysis. Given an algorithm whose iterates converge to a point $x_*$, extrapolation methods consist in combining iterates of the converging algorithm to obtain a new point hopefully closer to the solution. This idea was first applied for accelerating convergence of sequences in $\reals$ by fitting a linear model on the iterates and using as extrapolated point, the fixed point of this model \citep{Aitk27,Shan55,brezinski2006}. Extrapolation techniques have been extended to linearly converging sequences of vector \citep{Ande65,Sidi86} with convergence guarantees on auto-regressive models. More recently, those methods, and in particular Anderson acceleration, have been knowing a regain of interest in the optimization community. This renewed interest started with the work of \cite{Scie16}, which applied these extrapolation schemes to optimization algorithms, and where a regularization technique was proposed for obtaining convergence guarantees beyond auto-regressive settings. 

These extrapolation methods were widely extended: to the stochastic setting \citep{Scie17}, to composite optimization problems \citep{Mass18,Mai19}, for splitting methods \citep{Poon19,Fu19}, and to accelerate momentum based methods \citep{Boll18}.

Obtaining convergence guarantees for Anderson acceleration-type extrapolation, outside of the simple auto-regressive case (which corresponds to quadratic programs), is still an open issue. In \citep{Scie16}, explicit convergence guarantees are given in a regime asymptotically close to the solution. \cite{Zhan18} provides a globally converging Anderson acceleration type algorithm; however, due to the general setting of the paper, no convergence rate are provided. Nonasymptotic convergence bounds are provided in \citep{Toth15,Li20}, but they involve the inverse of the smallest eigenvalue of a Krylov matrix, usually very poorly conditioned. 

\textbf{Contribution :} Our contribution is two-fold.
\begin{itemize}
    \item We provide an explicit upper bound for the optimal value of the constrained Chebychev problem on polynomials. We demonstrate it is exact on some range of parameters and show numerically that it is close to the ground truth elsewhere.
    \item We us this to give an explicit linear convergence rate for constrained Anderson acceleration (CAA) scheme outside of the auto-regressive setting.
\end{itemize}
\subsection{Notations}
Depending on the context $\|\cdot\|$ either denotes the classical Euclidean norm, when applied to a vector in $\reals^n$, or the operator norm when applied to a matrix in $\reals^{n\times n}$. 

$\|\cdot\|_1$ is the sum of the absolute values of the components of a vector. When applied to a polynomial it is the sum of the absolute value of its coefficients.
\subsection{Setting}\label{sec:setting}
In this paper we study linear convergence of the constrained Anderson extrapolation scheme on some operator $F\;:\;\reals^n \rightarrow \reals^n$. 
For us, $F$ is typically a gradient step with constant stepsize. We make the following assumptions on $F$ throughout the paper.

\textbf{Assumptions:}
\begin{enumerate}
    \item $F$ is $\rho$ Lipschitz with $\rho < 1$. This implies that $F$ has a unique fixed point $x_*$ and the iterates of the fixed point iterations $x_{k+1} = F(x_k)$ converge with a linear rate $\rho$. 
    \item There is $G \in S_+^{n \times n}$ with $G \preccurlyeq \rho I$ and $\xi : \reals^n \rightarrow \reals^n$ $\alpha$-Lipschitz with $\alpha \geq 0$ such that $F = G + \xi$.
\end{enumerate}

In the case of $F$ encoding a gradient step, and given a $\mu$ strongly convex function $f : \reals^n \rightarrow \reals$, with $L$-Lipschitz gradient, it is well known that the operator $F = I - \tfrac{1}{L}\nabla f$ is contractive with $\rho = \left(1-\tfrac{\mu}{L} \right)$. In addition, assuming $f$ is $\mathcal{C}^2$ with $\eta$-Lipschitz hessian, we can show that $(I-\tfrac{1}{L}\nabla f)$ satisfies Assumption 2 around a point $x_0\in\reals^n$, with $\alpha$ proportional to $\eta\|\nabla f(x_0)\|$. This fact is made more precise in \S\ref{sec:main-res}.

We focus on the online version of Anderson acceleration \citep{Scie16,Scie18}. This means that the number of iterates used to perform the extrapolation is fixed to $k+1$ with some integer $k > 0$. 

\begin{algorithm}[!ht]
\caption{Constrained Anderson Acceleration}
\label{algo:CAA}
\begin{algorithmic}
\STATE \algorithmicrequire\;$x_0 \in \reals^n$, $F$ satisfying the assumptions, $C$ a constraint bound and $k$ controlling the number of extrapolation steps.
\FOR{$i = 0\dots k$}
\STATE $x_{i+1} = F(x_i)$
\ENDFOR
\STATE $R = \begin{bmatrix} x_0 -x_{1}&\cdots&x_{k}-x_{k+1}
\end{bmatrix}$ 
\STATE Compute \boxed{\tilde{c} = \underset{\substack{\mathbf{1}^Tc = 1,\; \|c\|_1 \leq C}}{\argmin}\;\|Rc\| } \hspace{\fill}(\refstepcounter{equation}{\theequation}\label{eq:CAA})
\STATE $x_e = \sum_{i=0}^{k}\tilde{c}_ix_{i}$
\STATE \algorithmicensure\; $x_e$
\end{algorithmic}
\end{algorithm}

Algorithm~\ref{algo:CAA} corresponds to one step of the online method. As we are interested in linear convergence rates, we look at results of the form $\|F(x_e)-x_e\| \leq \tilde{\rho}\|F(x_0)-x_0\|$, where $x_e$ is the output of Algorithm~\ref{algo:CAA} started at $x_0$. The quantity $\|F(x)-x\|$ is a criterion that is used to control how far is $x$ from being the fixed point of $F$, indeed $\|F(x)-x\| = 0 \iff x = x_*$. In addition, since $F$ is $\rho$ Lipschitz by assumption we also have that $\|F(x_{k})-x_{k}\| \leq \rho^{k}\|F(x_0)-x_0\|$ thus we say that extrapolation provides acceleration as soon as~$\tilde{\rho} < \rho^{k}$. 

It contrasts with the acceleration involving a convergence rate proportional to $(1-\sqrt{1-\rho})$ instead of $\rho$ (e.g Optimal Method of \cite{Nest18}), this type of acceleration is obtained for Anderson acceleration in the offline setting where $k$ grows to infinity, meaning that we use more and more iterates to perform extrapolation which is not recommended practice. 
\section{Constrained Anderson Acceleration}
We first recall some standard results on Anderson Acceleration applied to a linear operator when $\alpha = 0$. Then we introduce constraints on the extrapolation coefficients in order to deal with a perturbation parameter $\alpha > 0$.
\subsection{Anderson Acceleration on Linear Problem}
Let consider first the simple setting where $\alpha = 0$ and Algorithm~\ref{algo:CAA} is used with $C = \infty$. We recall the well-known Anderson acceleration result on the linear case.

\begin{proposition}
Let $F$ satisfying the assumptions of \S\ref{sec:setting} with $\alpha = 0$ and $x_e \in \reals^n$ the output of Algorithm~\ref{algo:CAA} started at $x_0\in\reals^n$ that is not the fixed point of $F$, with $C = \infty$ and $k+1 > 1$. We have that
\BEQS
\frac{\|F(x_e) - x_e\|}{\|F(x_0)-x_0\|} \leq \underset{\substack{p\in\reals_{k}[X]\\p(1)=1}}{\min}\; \underset{x\in [0,\rho]}{\max}|p(x)| = \rho_* := \tfrac{2\beta^k}{1+\beta^{2k}},
\EEQS
with $\beta = \tfrac{1-\sqrt{1-\rho}}{1+\sqrt{1-\rho}}$. In addition $\rho_* < \rho^{k}$.
\end{proposition}
\begin{proof}
Reformulation of \citep{Golu61,Scie16}
\end{proof}

\vspace{-0.2cm}
In the following, $\alpha$ may be nonzero and the previous Proposition does not apply.
\subsection{The Non Linear Case}
When applying the extrapolation step \eqref{eq:CAA} to a nonlinear operator $F$, the regularity of the matrix $R^TR \in \symm_{(k+1)}$ becomes an important issue. This matrix can be arbitrarily close to singular or even singular, for instance when $F$ is a gradient steps operator, it is known that the consecutive gradient tends to get aligned. Thus the solution vector $\tilde{c}$ can have coefficients with very large magnitude. When $\alpha>0$, those coefficients are multiplied with the nonlinear part of $F$ and can make the iterates of the algorithm diverge (see \citep{Scie16} for an example of such divergence). To fix this, one needs to control the magnitude of these coefficients, by e.g. regularizing problem~\eqref{eq:CAA}, as in \cite{Scie16} with $C=\infty$, or by imposing hard constraints on $\tilde{c}$, as we do. Regularizing involves easier computations in practice but imposing constraints makes the analysis simpler.

\begin{proposition}\label{prop:nonlinear-aa}
Given $F$ satisfying assumptions of \S\ref{sec:setting}, $\alpha \geq 0$ and $x_e \in \reals^n$ the output of Algorithm~\ref{algo:CAA} started at $x_0\in\reals^n$ with $C \geq 1$ and $k \geq 1$. We have
\BEQS
\|F(x_e)-x_e\| \leq \left(\max_{x\in[0,\rho]}|p_*(x)| + 3 C \alpha k \right)\|F(x_0)-x_0\|
\EEQS
where 
\[
p_* \in \underset{\substack{p \in \reals_k[X]\\ p(1)=1\\\|p\|_1\leq C}}{\argmin}\; \underset{x\in [0,\rho]}{\max} |p(x)|
\]
and $\|p\|_1$ is the $l_1$ norm of the vector of coefficients of~$p$.
\end{proposition}
\begin{proof} We provide here a sketch of the proof. A complete version is provided in Appendix~\ref{app:proof-nonlin-aa}.

Using the definition  $x_{i+1}= F(x_i) = G(x_i) + \xi(x_i)$.
one can show the following bound
\begin{align*}
    \lnorm F(x_e) - x_e\rnorm  \leq& \lnorm\sum_{i=0}^k\tilde{c}_i(x_{i+1}-x_i)\rnorm \\
    &+ \lnorm \xi(x_e) - \sum_{i=0}^k\tilde{c}_i\xi(x_i)\rnorm.
\end{align*}  
The first term of the right hand side is the quantity that is minimized in \eqref{eq:CAA} and the second term is due to the non linear part of $F$.

To control the first term we use the fact that since $\tilde{c}$ is solution of \eqref{eq:CAA} and $c^*$ the vector of coefficients of $p_*$ is admissible, then 
\[ \lnorm\sum_{i=0}^k\tilde{c}_i(x_{i+1}-x_i)\rnorm \leq \lnorm\sum_{i=0}^kc^*_i(x_{i+1}-x_i)\rnorm. \]

After some transformations using $\alpha$-Lischitzness of $\xi$, we obtain the bound wanted.
\end{proof}

\vspace{-0.3cm}
Proposition~\ref{prop:nonlinear-aa} exhibits a trade-off between (i) allowing coefficients to have larger magnitudes, i.e., via a large $C$, leading to a smaller $\underset{x\in[0,\rho]}{\max}|p_*(x)|$ that gets closer to $\rho_*$, and (ii) diminishing $C$ to better control the nonlinear part of $F$ but gettinga a slower rate  $\underset{x\in[0,\rho]}{\max}|p_*(x)|$ closer to~$\rho^k$. 

The following corollary simply states that one can allow a small relative error in the computation of \eqref{eq:CAA} and keep a linear convergence.
\begin{corollary}\label{cor:computing}
Under the conditions of Proposition~\ref{prop:nonlinear-aa}, if \eqref{eq:CAA} is solved with relative precision $\varepsilon\|F(x_0)-x_0\|$ on the optimal value for $\varepsilon > 0$ then 
\[\|(F-I)x_e\| \leq \left(\max_{x\in[0,\rho]}|p_*(x)| + 3 C \alpha k +\varepsilon \right)\|(F-I)x_0\|.\]
\end{corollary}

In the next section, we describe the variation of $\underset{x\in[0,\rho]}{\max}|p_*(x)|$ as a function of $C$, rendering the above trade-off explicit.

\section{Constrained Chebychev Problem}

\def \pltwidth = {.34}
\def \pltheight = {.305}
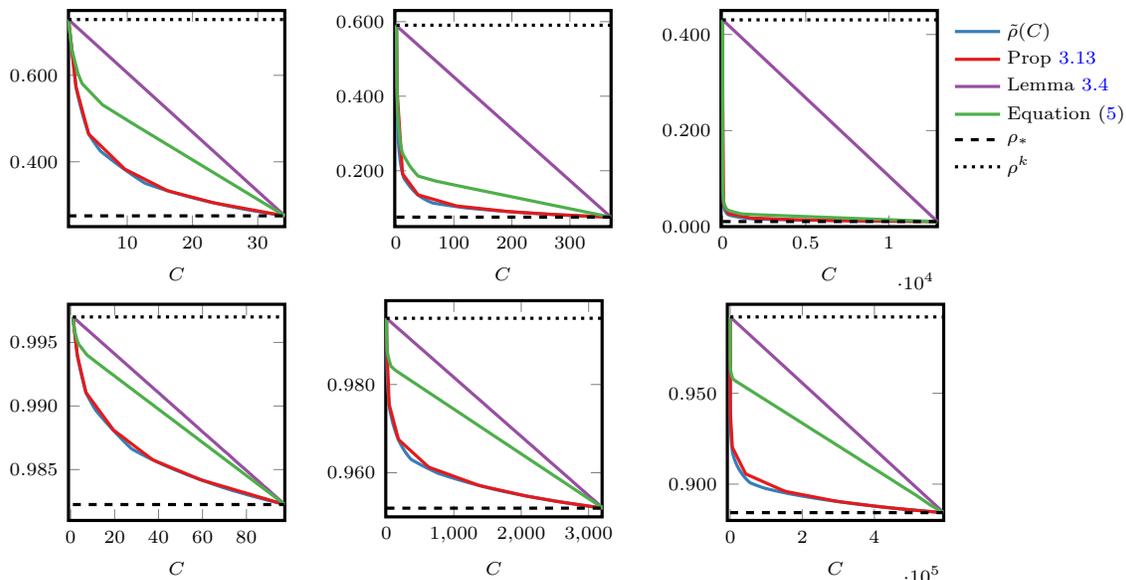
\begin{figure*}[!ht]
    \centering
    \begin{tabular}{ccc}
        \begin{tikzpicture}
    \begin{axis}[legend pos=outer north east, legend style={draw=none},legend cell align={left},xlabel={$C$}, plotOptions, ymin=0.25, ymax=0.75,xmin=1,xmax=34.141,width=0.26\linewidth, height=0.26\linewidth,y tick label style={
        /pgf/number format/.cd,
        fixed,
        fixed zerofill,
        precision=3,
        /tikz/.cd
    }]
		\addplot [color=colorP1] table [x=C,y=ex] {Fig/0.9_3.txt};
		\addplot [color=colorP2] table [x=C,y=up] {Fig/0.9_3.txt};
		\addplot [color=colorP3, domain=1:34.141, samples=50]
		{(x-1)/(34.141-1)*2*((1-sqrt(0.1))/(1+sqrt(0.1)))^3/(1+((1-sqrt(0.1))/(1+sqrt(0.1)))^6) + (34.141-x)/(34.141-1)*0.9^3};
		\addplot [color=colorP4] table [x=C,y=proj] {Fig/0.9_3.txt};
		\addplot [color=black,dashed ,domain=1:34.141, samples=50] {2*((1-sqrt(0.1))/(1+sqrt(0.1)))^3/(1+((1-sqrt(0.1))/(1+sqrt(0.1)))^6)};
		\addplot [color=black,dotted,domain=1:34.141, samples=50] {0.9^3};
		\end{axis}
		
	\end{tikzpicture}
    		&
		\begin{tikzpicture}
		\begin{axis}[legend pos=outer north east, legend style={draw=none},legend cell align={left},xlabel={$C$}, plotOptions, ymin=0.05, ymax=0.63,xmin=-2,xmax=371,width=0.26\linewidth, height=0.26\linewidth,y tick label style={
        /pgf/number format/.cd,
        fixed,
        fixed zerofill,
        precision=3,
        /tikz/.cd
    }]
		\addplot [color=colorP1] table [x=C,y=ex] {Fig/0.9_5.txt};
		\addplot [color=colorP2] table [x=C,y=up] {Fig/0.9_5.txt};
		\addplot [color=colorP3, domain=1:370.62, samples=50] {(x-1)/(370.62-1)*2*((1-sqrt(0.1))/(1+sqrt(0.1)))^5/(1+((1-sqrt(0.1))/(1+sqrt(0.1)))^10) + (370.62-x)/(370.62)*0.9^5};
		\addplot [color=colorP4] table [x=C,y=proj] {Fig/0.9_5.txt};
		\addplot [color=black,dashed ,domain=1:370.62, samples=50] {2*((1-sqrt(0.1))/(1+sqrt(0.1)))^5/(1+((1-sqrt(0.1))/(1+sqrt(0.1)))^10)};
		\addplot [color=black,dotted ,domain=1:370.62, samples=50] {0.9^5};
		\end{axis}
	\end{tikzpicture}     & 
	\begin{tikzpicture}
    \begin{axis}[legend pos=outer north east, legend style={draw=none},legend cell align={left},xlabel={$C$}, plotOptions, ymin=0, ymax=0.45,xmin=-100,xmax=12920,width=0.26\linewidth, height=0.26\linewidth,y tick label style={
        /pgf/number format/.cd,
        fixed,
        fixed zerofill,
        precision=3,
        /tikz/.cd
    }]
		\addplot [color=colorP1] table [x=C,y=ex] {Fig/0.9_8.txt};
		\addplot [color=colorP2] table [x=C,y=up] {Fig/0.9_8.txt};
		\addplot [color=colorP3, domain=1:12920, samples=50] {(x-1)/(12920-1)*2*((1-sqrt(0.1))/(1+sqrt(0.1)))^8/(1+((1-sqrt(0.1))/(1+sqrt(0.1)))^16) + (12920-x)/(12920)*0.9^8};
		\addplot [color=colorP4] table [x=C,y=proj] {Fig/0.9_8.txt};
		\addplot [color=black,dashed ,domain=1:12920, samples=50] {2*((1-sqrt(0.1))/(1+sqrt(0.1)))^8/(1+((1-sqrt(0.1))/(1+sqrt(0.1)))^16)};
		\addplot [color=black,dotted ,domain=1:12920, samples=50] {0.9^8};
		\addlegendentry{$\tilde{\rho}(C)$}
		\addlegendentry{Prop~\ref{prop:global-bound}}
		\addlegendentry{Lemma~\ref{lem:cvx-cheb}}
		\addlegendentry{Equation~\eqref{eq:proj-inf-norm}}
		\addlegendentry{$\rho_*$}
		\addlegendentry{$\rho^k$}
		\end{axis}
		
	\end{tikzpicture} 
	\\
	\begin{tikzpicture}
		\begin{axis}[legend pos=outer north east, legend style={draw=none},legend cell align={left},xlabel={$C$}, plotOptions, ymin=0.981, ymax=0.998,xmin=-1,xmax=97.45,width=0.26\linewidth, height=0.26\linewidth,y tick label style={
        /pgf/number format/.cd,
        fixed,
        fixed zerofill,
        precision=3,
        /tikz/.cd
    }]
		\addplot [color=colorP1] table [x=C,y=ex] {Fig/0.999_3.txt};
		\addplot [color=colorP2] table [x=C,y=up] {Fig/0.999_3.txt};
		\addplot [color=colorP3, domain=1:97.45, samples=50] {(x-1)/(97.45-1)*2*((1-sqrt(0.001))/(1+sqrt(0.001)))^3/(1+((1-sqrt(0.001))/(1+sqrt(0.001)))^6) + (97.45-x)/(97.45-1)*0.999^3};
		\addplot [color=colorP4] table [x=C,y=proj] {Fig/0.999_3.txt};
		\addplot [color=black,dashed ,domain=1:97.45, samples=50] {2*((1-sqrt(0.001))/(1+sqrt(0.001)))^3/(1+((1-sqrt(0.001))/(1+sqrt(0.001)))^6)};
		\addplot [color=black,dotted ,domain=1:97.45, samples=50] {0.999^3};
		\end{axis}
	\end{tikzpicture}
	&
        \begin{tikzpicture}
		\begin{axis}[legend pos=outer north east, legend style={draw=none},legend cell align={left},xlabel={$C$}, plotOptions, ymin=0.95, ymax=0.999,xmin=-10,xmax=3200,width=0.26\linewidth, height=0.26\linewidth,y tick label style={
        /pgf/number format/.cd,
        fixed,
        fixed zerofill,
        precision=3,
        /tikz/.cd
    }]
		\addplot [color=colorP1] table [x=C,y=ex] {Fig/0_999_5.txt};
		\addplot [color=colorP2] table [x=C,y=up] {Fig/0_999_5.txt};
		\addplot [color=colorP3, domain=1:3213, samples=50] {(x-1)/(3213.8-1)*2*((1-sqrt(0.001))/(1+sqrt(0.001)))^5/(1+((1-sqrt(0.001))/(1+sqrt(0.001)))^10) + (3213.8-x)/(3213.8-1)*0.999^5};
		\addplot [color=colorP4] table [x=C,y=proj] {Fig/0.999_5.txt};
		\addplot [color=black,dashed ,domain=1:3213, samples=50] {2*((1-sqrt(0.001))/(1+sqrt(0.001)))^5/(1+((1-sqrt(0.001))/(1+sqrt(0.001)))^10)};
		\addplot [color=black,dotted ,domain=1:3213, samples=50] {0.999^5};
		\end{axis}
	\end{tikzpicture} &
	\hspace{-2.6cm}\vspace{-0.5cm}
	\begin{tikzpicture}
		\begin{axis}[legend pos=outer north east, legend style={draw=none},legend cell align={left},xlabel={$C$}, plotOptions, ymin=0.88, ymax=0.999,xmin=-8000,xmax=5.9216e+05,width=0.26\linewidth, height=0.26\linewidth,y tick label style={
        /pgf/number format/.cd,
        fixed,
        fixed zerofill,
        precision=3,
        /tikz/.cd
    }]
		\addplot [color=colorP1] table [x=C,y=ex] {Fig/0.999_8.txt};
		\addplot [color=colorP2] table [x=C,y=up] {Fig/0.999_8.txt};
		\addplot [color=colorP3, domain=1:5.9216e+05, samples=50] {(x-1)/(5.9216e+05-1)*2*((1-sqrt(0.001))/(1+sqrt(0.001)))^8/(1+((1-sqrt(0.001))/(1+sqrt(0.001)))^16) + (5.9216e+05-x)/(5.9216e+05-1)*0.999^8};
		\addplot [color=colorP4] table [x=C,y=proj] {Fig/0.999_8.txt};
		\addplot [color=black,dashed ,domain=1:5.9216e+05, samples=50] {2*((1-sqrt(0.001))/(1+sqrt(0.001)))^8/(1+((1-sqrt(0.001))/(1+sqrt(0.001)))^16)};
		\addplot [color=black,dotted ,domain=1:5.9216e+05, samples=50] {0.999^8};
		\end{axis}
	\end{tikzpicture}
    \end{tabular}
    \caption{Blue curves correspond to numerical solutions to \eqref{eq:sos}, red ones correspond to the bound from Proposition~\ref{prop:global-bound} with $M=k$, purple ones are the bound from Lemma~\ref{lem:cvx-cheb} and the greens are the maximum of the absolute value of the solution of \eqref{eq:proj-inf-norm} over $[0,\rho]$. On x-axis $C$ goes from $1$ to $C_*$ defined in \eqref{eq:c*}. Top : $\rho=0.9$. Bottom: $\rho=0.999$. Left: $k=3$. Middle: $k=5$. Right: $k=8$.}
  \label{fig:sos-lin}
\end{figure*}

The Chebyshev problem, defined in the following theorem, is fundamental in many field of numerical analysis. It is used to provide convergence rate for many optimization methods such as the conjugate gradient algorithm, Anderson acceleration, or Chebyshev iterations \citep{Golu61,Nemi84a,Shew94,Nemi92}. 
\begin{theorem}[\cite{Golu61}]\label{thm:recall-cheb}
Let $\rho > 0$ and $k > 0$, we call Chebyshev problem of degree $k$ on $[0,\rho]$ the following optimization problem
\BEQ\label{eq:non-const-cheb}
\underset{\substack{p \in \reals_k[X]\\p(1)=1}}{\min\;}\underset{x\in[0,\rho]}{\max\;}|p(x)| \tag{Cheb}
\EEQ
The solution of this optimization problem is $p_*(X) = \tfrac{T_k(\tfrac{2X-\rho}{\rho})}{\left|T_k(\tfrac{2-\rho}{\rho})\right|}$ where $T_k$ is the first kind Chebyshev polynomial of order $k$. The optimal value is equal to $\rho_*$ defined in Proposition~\ref{prop:nonlinear-aa}.
\end{theorem}
\begin{proof}
For completeness, a proof of this result is provided in Appendix~\ref{prop:rescaled-cheb}. 
\end{proof}

\vspace{-0.2cm}
We have seen in Proposition~\ref{prop:nonlinear-aa} that we need to control the optimal value of a slightly modified problem where we add a constraint on the $l_1$ norm of the vector of coefficients of the polynomial. Adding this constraint breaks the explicit result of Theorem~\ref{thm:recall-cheb}. No closed form solution for the constrained Chebyshev problem is known in the general case. In this section we aim at providing good upper bounds on the optimal value of this constrained problem.

Given $k \geq 1$ and $\rho \in ]0,1[$, we are now interested in the following constrained Chebyshev problem 
\BEQ\label{eq:ctr-cheb}
\tilde{\rho}(C) := \underset{\substack{p \in \reals_k[X]\\ p(1)=1\\\|p\|_1\leq C}}{\min\;} \underset{x\in [0,\rho]}{\max} |p(x)| \tag{Ctr-Cheb}
\EEQ
The next subsection describes how to compute this $\tilde{\rho}(C)$ numerically for $C \geq 1$. Note that the admissible set is empty when $C < 1$.

\subsection{Numerical Solutions}
When $C \geq 1$, the problem \eqref{eq:ctr-cheb} has a non empty admissible set. In addition, the feasible set an intersection of an affine space and an $l_1$ ball, and hence is convex. The objective function being a norm on $\reals_k[X]$ and equivalently on $\reals^{k+1}$, is convex.
The problem \eqref{eq:ctr-cheb} is equivalent to 
\BEQ
\begin{array}{ll}
     &\mbox{min}\quad t  \\
     &p \in \reals_k[X],\;t\in \reals\\
     &p(1)=1,\;\|p\|_1\leq C\\
     &-t \leq p(x) \leq t,\;\forall x \in [0,\rho]\\
     
\end{array} 
\EEQ
This problem involves positivity constraints of polynomials on a bounded interval. A classical argument to transform this local positivity into positivity on $\reals$ is using the following change of variable. $p(x) \geq 0\; \forall x\in[0,\rho] \iff (1+x^2)^k p\left(\rho\tfrac{x^2}{1+x^2}\right) \geq 0 \; \forall x \in \reals$. Then, positiviy constraint for a polynomial on $\reals$ can be relaxed using a sum of squares (SOS) formulation~\citep{parrilo2000,lasserre01}, which is exact for univariate polynomials (e.g \citep[Theorem 1]{Magr19} for a short proof). Standard packages can be used to solve efficiently the following optimization problem with SOS constraints. 
\BEQ\label{eq:sos}
\begin{array}{ll}
     &\mbox{min}\quad t  \\
     &p \in \reals_k[X],\;t\in \reals\\
     &p(1)=1,\;\|p\|_1\leq C\\
     &(1+x^2)^kp(\rho\tfrac{x^2}{1+x^2}) + (1+x^2)^kt \geq 0\;\forall x \in \reals\\
     &(1+x^2)^kt-(1+x^2)^kp(\rho\tfrac{x^2}{1+x^2})  \geq 0\;\forall x \in \reals
     
\end{array} 
\EEQ
We used YALMIP~\citep{Lof04} and MOSEK~\citep{mosek}. Numerical solutions to~\eqref{eq:sos} are provided in Figure~\ref{fig:sos-lin} (blue) for a few values of $\rho$ and $k$.

\subsection{Exact Bounds and Upper Bounds}
The main goal of this part is to provide an explicit upper bound for the function $\tilde{\rho}(C)$ defined in \eqref{eq:ctr-cheb}, for using it in the result of Proposition~\ref{prop:nonlinear-aa}. 

We denote by $C_*$ the $l_1$-norm of the rescaled Chebyshev polynomial $p_*$ of Theorem~\ref{thm:recall-cheb}. 
\BEQ\label{eq:c*} 
C_* = \|p_*\|_1,\text{ where } p_* \text{ solves } \eqref{eq:non-const-cheb}
\EEQ
\begin{remark}\label{rem:cst}
By Theorem~\ref{thm:recall-cheb}, $\tilde{\rho}(C)$ is constant equal to $\rho_*$ as soon as $C \geq C_*$.
\end{remark}
\begin{remark}
When $C = 1$, the admissible set of \eqref{eq:ctr-cheb} consists only on the monomials of degree smaller than $k$. $X^k$ has the minimal absolute value on $[0,\rho]$ and $\tilde{\rho}(1) = \rho^k$.
\end{remark}

The following lemma gives a first upper bound on $\tilde{\rho}$.

\begin{lemma}\label{lem:cvx-cheb}
The function $\tilde{\rho}$ defined in \eqref{eq:ctr-cheb} is convex on $[1,+\infty[$.
Thus for $C \in [1,C_*]$ 
\[ \tilde{\rho}(C) \leq \tfrac{C_*-C}{C_*-1}\rho^k + \tfrac{C-1}{C_*-1}\rho_* \]
\end{lemma}

This is actually a very coarse bound due to the fact that $C_* >> 1$. Indeed we can observe in Figure~\ref{fig:sos-lin} that there is a important gap between $\tilde{\rho}$ and the coarse upper bound from Lemma~\ref{lem:cvx-cheb} that is displayed in purple. However we will use the convexity of $\tilde{\rho}$ to link together finer upper bounds taken at different $C$'s.

Another upper bound is represented on Figure~\ref{fig:sos-lin} (green). This is the maximum over $[0,\rho]$ of the absolute value of $q$ defined as
\BEQ\label{eq:proj-inf-norm}
q \in \underset{\substack{p(1)=1\\ \|p\|_1\leq C}}{\argmin}\; \underset{x\in[0,\rho]}{\max}\left|p(x) - p_*(x)\right|
\EEQ
where $p_*$ is solution of \eqref{eq:non-const-cheb}. This is a naive upper bound but it is as hard to compute as solving \eqref{eq:ctr-cheb}.

In the next lemma, we provide an explicit expression of $\tilde{\rho}(C)$ for $C$ in an explicit neighbourhood of $1$.
\begin{lemma}\label{lem:tilde-rho-1} For $C \in [1,\tfrac{2+\rho^k}{2-\rho^k}]$ we have the following expression for $\tilde{\rho}$
\[\tilde{\rho}(C) = \tfrac{C+1}{2}\rho^k - \tfrac{C-1}{2}  \]
\end{lemma}
\begin{proof}
Let show that $p(X) = \tfrac{C+1}{2}X^k - \tfrac{C-1}{2}$ is solution of \eqref{eq:ctr-cheb}. First notice that $p$ is feasible as $\|p\|_1 = C$ and $p(1) = 1$. In addition, since $C \in [1,\tfrac{2+\rho^k}{2-\rho^k}]$ we have $\underset{x\in[0,\rho]}{\max}\; |p(x)| = p(\rho) = \tfrac{C+1}{2}\rho^k - \tfrac{C-1}{2}$.

Let $q$ be another feasible polynomial such that $q = \sum_{i=0}^kq_iX^i$ and $q \neq p$. We show that $|q(\rho)|\geq |p(\rho)| $.
\begin{align*}
    q(\rho) &= \sum_{c_i \geq 0}q_i\rho^i + \sum_{c_i \leq 0}q_i\rho^i\\
    & \geq \sum_{c_i \geq 0}q_i\rho^k +  \sum_{c_i \leq 0}q_i\\
    & = \sum_{c_i \geq 0}q_i\rho^k +  \left(1-\sum_{c_i \geq 0}q_i\right)
\end{align*}
In addition one notices that $\sum_{c_i \geq 0}q_i - \sum_{c_i \leq 0}q_i \leq C$ thus $\sum_{c_i \geq 0}q_i \leq  \tfrac{C+1}{2}$ and thus
\begin{align*}
    q(\rho) & \geq \tfrac{C+1}{2}(\rho^k-1) + 1\\
    &= p(\rho) > 0
\end{align*}
Thus $|q(\rho)| = q(\rho) \geq p(\rho) = \underset{x\in[0,\rho]}{\max}\; |p(x)|$. Then $\underset{x\in[0,\rho]}{\max}\; |q(x)| \geq \underset{x\in[0,\rho]}{\max}\; |p(x)|$ and necessarily $p$ is a solution of \eqref{eq:ctr-cheb}.
\end{proof}

\begin{remark}
When $k=1$, $C_* = \tfrac{2+\rho}{2-\rho}$, and the function $\tilde{\rho}$ is entirely defined by Remark~\ref{rem:cst} and Lemma~\ref{lem:tilde-rho-1}.
\end{remark}
\begin{remark}
We have $\tilde{\rho}(\tfrac{2+\rho^k}{2-\rho^k}) = \tfrac{\rho^k}{2-\rho^k}$. When $\rho \rightarrow 1$, $\tfrac{2+\rho^k}{2-\rho^k} \rightarrow 3$, and  $\tfrac{1-\tfrac{\rho^k}{2-\rho^k}}{1-\rho^k} \rightarrow 2$. This means $\tilde{\rho}(\tfrac{2+\rho^k}{2-\rho^k}) = 1 - 2(1-\rho^k) + o(1-\rho^k)$.
\end{remark}

The following proposition gives the form of the solution of \eqref{eq:ctr-cheb} in a neighborhood of  $C_*$. It states that the solutions of the Chebyshev problem with light constraints are also rescaled Chebyshev polynomials on a segment $[-\varepsilon,\rho]$ instead of $[0,\rho]$, with $\varepsilon\geq0$.

\begin{proposition} \label{prop:cheb-pb-eps} Let $\rho\in]0,1[$, and $k \geq 1$. For any $\varepsilon \in [0,\tilde{\varepsilon}]$ with $\tilde{\varepsilon} = \rho\tfrac{1+\cos(\tfrac{2k-1}{2k}\pi)}{1-\cos(\tfrac{2k-1}{2k}\pi)}$ we have \[\tilde{\rho}(\|p_\varepsilon\|_1) = \underset{x\in[-\varepsilon,\rho]}{\max}\; |p_\varepsilon(x)| \] 
where $p_\varepsilon = \underset{\substack{p\in\reals_k[X]\\p(1)=1}}{\argmin}\;\underset{x\in[-\varepsilon,\rho]}{\max}\; |p(x)|  $
\end{proposition}
\begin{proof} The proof is provided in Appendix~\ref{app:cheb-pb-eps}.
\end{proof}


\vspace{-0.2cm}
Due to the fact that the coefficients of the polynomials $p_\varepsilon$ for $\varepsilon\in[0,\tilde{\varepsilon}]$ have alternating signs, we can get a nice expression for $\|p_\varepsilon\|_1$.
\begin{lemma}\label{lem:norm-cheb-eps}  Let $\rho\in]0,1[$, and $k \geq 1$. For any $\varepsilon \in [0,\tilde{\varepsilon}]$ with $\tilde{\varepsilon} = \rho\tfrac{1+\cos(\tfrac{2k-1}{2k}\pi)}{1-\cos(\tfrac{2k-1}{2k}\pi)} $ we have
\begin{align*}
    \|p_\varepsilon\|_1  &= \tfrac{\tilde{\rho}_\varepsilon}{2}\left(\tfrac{1}{\rho+\varepsilon}\right)^k\left(\left(2 +\rho-\varepsilon -2\sqrt{(1+\rho)(1-\varepsilon)}\right)^k \right.\\
    & \left.+\left(2 +\rho-\varepsilon +2\sqrt{(1+\rho)(1-\varepsilon)}\right)^k \right)
\end{align*}
where $\tilde{\rho}_\varepsilon = \tfrac{2\beta_\varepsilon^k}{1+\beta_\varepsilon^{2k}}$ and $\beta_\varepsilon = \tfrac{1-\sqrt{1-\tfrac{\rho+\varepsilon}{1+\varepsilon}}}{1+\sqrt{1-\tfrac{\rho+\varepsilon}{1+\varepsilon}}}$.

\end{lemma}
\begin{proof}
From the proof in Appendix~\ref{app:cheb-pb-eps} we can see that the coefficients of $p_\varepsilon$ alternate signs, which means that $\|p_\varepsilon\|_1 = |p_\varepsilon(-1)|$. We then use the classical expression for $T_k(x)$ with $|x| \geq 1$ (see e.g \citep[Eq 1.49]{Maso02})\\
$T_k(x) = \tfrac{1}{2}\left((x-\sqrt{x^2-1})^k + (x+\sqrt{x^2-1})^k \right) .$
\end{proof}
\vspace{-0.3cm}
\begin{remark}\label{rem:c*} In particular, one can recover the value of $C_\star$ (which we did not provide before) from Lemma~\ref{lem:norm-cheb-eps} applied to the unconstrained Chebyshev problem~\eqref{eq:ctr-cheb} ($\varepsilon =0$)
From previous lemma we observe that 
\[C_* = \tfrac{\rho_*}{2\rho^k}\left((2+\rho -2\sqrt{1+\rho})^k + (2+\rho +2\sqrt{1+\rho})^k \right) \]
\end{remark}

Proposition~\ref{prop:cheb-pb-eps} and Lemma~\ref{lem:norm-cheb-eps} do not provide direct access to $\tilde{\rho}(C)$ for $C\in [\tilde{C},C_*]$ with $\tilde{C}=\|p_{\tilde{\varepsilon}}\|_1$. Indeed, we cannot explicitly invert the relation $\varepsilon \rightarrow \|p_\varepsilon\|_1$. However one can get arbitrarily good upper bounds by sampling $(\varepsilon_i)_{i\in[1,M]} \in [0,\tilde{\varepsilon}]$. Then, one can compute $C_i = \|p_{\varepsilon_i}\|_1$ explicitly using Lemma~\ref{lem:norm-cheb-eps}, and use convexity to interpolate linearly between the $C_i$ and $\tilde{\rho}(C_i)$, reaching a piecewise linear upper bound on $[\tilde{C},C]$. Arbitrarily high accuracy can be obtained by increasing $M$ in this procedure.

To construct the following upper bound to $\tilde{\rho}(C)$ on the remaining $C$'s, we rely on some numerical observations. Indeed we observed that the maximum over $[0,\rho]$ of rescaled Chebyshev polynomials $p_\varepsilon$ provide a good upper bound of $\tilde{\rho}(\|p_\varepsilon\|_1)$ for a large range of $\varepsilon$. In particular, we study $p_\rho$ in the following as we can get a relatively simple expression for $\|p_\rho\|_1$.

\begin{proposition}
Let $\rho \in ]0,1[$ and $k \geq 1$,
\[
\tilde{\rho}(C_1) \leq \rho_1 \triangleq \tfrac{2\beta_\rho^k}{1+\beta_\rho^{2k}}  
\]
with
\BEQ\label{eq:rho1}
C_1=\tfrac{\rho_1}{2\rho^k}\left((1-\sqrt{1+\rho^2})^k + (1+\sqrt{1+\rho^2})^k \right),
\EEQ and $\beta_\rho =
\tfrac{\sqrt{1+\rho}-\sqrt{1-\rho}}{\sqrt{1+\rho}+\sqrt{1-\rho}}$.
\end{proposition}
\begin{proof}
By Lemma~\ref{lem:rescaled-cheb-eps}, $p_\rho(X) = \rho_1T_k(\tfrac{X}{\rho})$. 
Let show that $\|p_\rho\|_1 = C_1$. To see that one can notice that $\|p_\varepsilon\|_1 = \rho_1\left|T_k(\tfrac{i}{\rho}) \right| = \tfrac{1}{2\rho^k}\left|\left(i - (i^2-\rho)^{1/2}\right)^k + \left(i - (i^2-\rho)^{1/2}\right)^k)\right| = C_1$, using \citep[Eq 1.49]{Maso02}. 

Thus $p_\rho$ is admissible for \eqref{eq:ctr-cheb} and its maximum absolute value on $[-\rho,\rho]$, $\rho_1$ provides an upper bound for $\tilde{\rho}(C_1)$.
\end{proof}

\vspace{-0.2cm}
We study the regime $\rho\sim 1$ in order to get more insight on what the previous results mean.
\begin{remark} 
We observe that 
\[\rho_1 = \tfrac{2\rho^k}{(1+\sqrt{1-\rho^2})^k+(1-\sqrt{1-\rho^2})^k}\leq \rho^k \]
and when $\rho \rightarrow 1$
\BEQS
\begin{aligned}
\rho_1-\rho_*&\sim  \tfrac{k}{2k-1}(\rho^k-\rho_*)\\
C_1 &\sim \tfrac{(1+\sqrt{2})^k + (1-\sqrt{2})^k}{(1+\sqrt{2})^{2k} + (1-\sqrt{2})^{2k}}C_* \leq \tfrac{1}{2^k}C_*    
\end{aligned}
\EEQS

In the regime $\rho \sim 1$ we have that by decreasing the constraint $C$ by a factor $2^k$, one only looses a factor $\tfrac{k}{2k-1}$ of possible acceleration.
\end{remark}

We propose the following upper bound for the function $\tilde{\rho}(C)$ over $[1,\infty[$.

\begin{proposition}\label{prop:global-bound}
Let $C \geq 1$, $k > 2$ , $0<\rho<1$. Choosing $M\in N^*$, $(\varepsilon_i)_{i\in[1,M]} = \left(\tfrac{\rho}{2^{i-1}}\right)_{i\in[1,M]}$ . Let $(C_i)_{i\in[1,M]} = (\|p_{\varepsilon_i}\|_1)_{i\in[1,M]}$, $C_{-1} = 1$, $C_0 = \tfrac{2+\rho^k}{2-\rho^k}$ and $C_{M+1}=C_*$. We denote by $(\rho_i)_{i\in[1,M]} = \left(\tfrac{2\beta_i^k}{1+\beta_i^{2k}} \right) $ with $\beta_i = \tfrac{1-\sqrt{1-\tfrac{\rho+\varepsilon_i}{1+\varepsilon_i}}}{1+\sqrt{1-\tfrac{\rho+\varepsilon_i}{1+\varepsilon_i}}}$, $\rho_{-1} = \rho^k$, $\rho_0 = \tfrac{\rho^k}{2-\rho^k}$ and $\rho_{M+1} = \rho_*$. Then
\[\tilde{\rho}(C) \leq \max_{i\in[-1,M]}\left(\tfrac{C-C_i}{C_{i+1}-C_i}\rho_{i+1} + \tfrac{C_{i+1}-C}{C_{i+1}-C_i}\rho_i,\rho_*\right) \]
is an explicit upper bound on $\tilde{\rho}$.
\end{proposition}
\begin{proof}
The $C_i = \|p_{\varepsilon_{i}}\|_1$ can be made explicit using for instance the formula in \citep[Equation 2.18]{Maso02} for the coefficients of $T_k$. If $C > C_*$ we saw that $\tilde{\rho}(C)=\rho_*$. Otherwise $C$ is between a $C_i$ and a $C_{i+1}$ and the result follows from the convexity of $\tilde{\rho}$.
\end{proof}

\vspace{-0.2cm}
As shown in Figure~\ref{fig:sos-lin}, Proposition~\ref{prop:global-bound} provides partially numerical upper bounds, represented in red on the figure, that are close to $\tilde{\rho}(C)$.  In the next section we use this result to provide explicit bounds on constrained Anderson acceleration for non linear operators.
\section{Main Results}\label{sec:main}
As discussed above, combining Proposition~\ref{prop:global-bound} and Proposition~\ref{prop:nonlinear-aa} gives an explicit linear rate of convergence for one pass of Algorithm~\ref{algo:CAA}.  
\begin{proposition}\label{prop:bound-aa-tot}
Let $F$ be satisfying assumptions of \S\ref{sec:setting}, $\alpha \geq 0$, and $x_e \in \reals^n$ be the output of Algorithm~\ref{algo:CAA} with $x_0\in\reals^n$, $C \geq 1$, and $k > 2$. It holds that
\[\|F(x_e)-x_e\| \leq \hat{\rho}(C)\|F(x_0)-x_0\| \]
where $\hat{\rho} = \tilde{\rho}(C) + 3\alpha k C$ with $\tilde{\rho}$ defined in \eqref{eq:ctr-cheb}.
In addition, given $M\in\N^*$,
\[\hat{\rho}(C) \leq \max_{i\in[-1,M]}\left(\tfrac{C-C_i}{C_{i+1}-C_i}\rho_{i+1} + \tfrac{C_{i+1}-C}{C_{i+1}-C_i}\rho_i,\rho_*\right) +3\alpha k C \]
where the $\rho_i$, $C_i$ are defined in Proposition~\ref{prop:global-bound}.
\end{proposition}

In what follows, we explicitly control this linear rate to provide guarantees on acceleration.

\subsection{Explicit Upper Bound on CAA Convergence Rate}
For readability purposes, let us consider a very simple form of the upper bound of $\tilde{\rho}(C)$, using Proposition~\ref{prop:global-bound} with $M=1$. This corresponds to a piecewise linear upper bound on $\tilde{\rho}(C)$ with four parts. The following proposition provides a range of values of $C$ for which Algorithm~\ref{algo:CAA} accelerates the convergence guarantees of the iterative process $F$, depending on the perturbation parameter~$\alpha$.

\begin{proposition}\label{prop:bound-aa}
Under the assumptions of Proposition~\ref{prop:bound-aa-tot},
\begin{enumerate}[label=(\roman*)]
    \item As soon as $\alpha < \tfrac{\rho^k(1-\rho^k)}{3k(2+\rho^k)} := \alpha_0$ there exists  a non empty interval $I$ containing $\tfrac{2+\rho^k}{2-\rho^k}$ such that $\hat{\rho}(C) < \rho^k$ for $C\in I$. 
    \item $\alpha < \min(\tfrac{\rho^k(1-\rho^k)}{3k(2+\rho^k)},\tfrac{\rho^k-\rho_1}{3kC_1}) := \alpha_1$ \\
    $\implies \hat{\rho}(C) < \rho^k \text{ for } C\in[\tfrac{2+\rho^k}{2-\rho^k},C_1] $.
    \item $\alpha < \min(\tfrac{\rho^k(1-\rho^k)}{3k(2+\rho^k)},\tfrac{\rho^k-\rho_*}{3kC_*}) := \alpha_2$\\ $\implies \hat{\rho}(C) < \rho^k \text{ for } C\in[\tfrac{2+\rho^k}{2-\rho^k},C_*]$.
\end{enumerate}
\end{proposition}
\vspace{-0.3cm}
\begin{proof}
See Appendix~\ref{app:bound-aa}
\end{proof}

\begin{figure*}[!ht]
    \centering
    \begin{tabular}{ccc}
    \begin{tikzpicture}
		\begin{semilogxaxis}[legend pos=outer north east, legend style={draw=none},legend cell align={left},xlabel={$C$}, plotOptions, ymin=0.948, ymax=0.999,xmin=1,xmax=10000,width=0.265\linewidth, height=0.275\linewidth]
		\addplot [color=colorP1] table [x=C,y=a0] {Fig/res_best0.999_5.txt};
		\addplot [color=colorP2] table [x=C,y=ae-9] {Fig/res_best0.999_5.txt};
		\addplot [color=colorP3] table [x=C,y=ae-8] {Fig/res_best0.999_5.txt};
		\addplot [color=colorP4] table [x=C,y=ae-7] {Fig/res_best0.999_5.txt};
		\addplot [color=colorP5] table [x=C,y=ae-6] {Fig/res_best0.999_5.txt};
		\addplot [color=colorP6] table [x=C,y=ae-5] {Fig/res_best0.999_5.txt};
		
		\addplot [color=black,dashed ,domain=1:10000, samples=50] {2*((1-sqrt(0.001))/(1+sqrt(0.001)))^5/(1+((1-sqrt(0.001))/(1+sqrt(0.001)))^10)};
		\addplot [color=black,dotted ,domain=1:10000, samples=50] {0.999^5};
		\end{semilogxaxis}
	\end{tikzpicture}     &
	\begin{tikzpicture}
		\begin{semilogxaxis}[legend pos=outer north east, legend style={draw=none},legend cell align={left},xlabel={$C$}, plotOptions, ymin=0.948, ymax=0.999,xmin=1,xmax=10000,width=0.265\linewidth, height=0.275\linewidth]
		\addplot [color=colorP1] table [x=C,y=a0] {Fig/res_up_with_k_parts0.999_5.txt};
		\addplot [color=colorP2] table [x=C,y=ae-9] {Fig/res_up_with_k_parts0.999_5.txt};
		\addplot [color=colorP3] table [x=C,y=ae-8] {Fig/res_up_with_k_parts0.999_5.txt};
		\addplot [color=colorP4] table [x=C,y=ae-7] {Fig/res_up_with_k_parts0.999_5.txt};
		\addplot [color=colorP5] table [x=C,y=ae-6] {Fig/res_up_with_k_parts0.999_5.txt};
		\addplot [color=colorP6] table [x=C,y=ae-5] {Fig/res_up_with_k_parts0.999_5.txt};
		
		\addplot [color=black,dashed ,domain=1:10000, samples=50] {2*((1-sqrt(0.001))/(1+sqrt(0.001)))^5/(1+((1-sqrt(0.001))/(1+sqrt(0.001)))^10)};
		\addplot [color=black,dotted ,domain=1:10000, samples=50] {0.999^5};
		\end{semilogxaxis}
	\end{tikzpicture} &
	\begin{tikzpicture}
		\begin{semilogxaxis}[legend pos=outer north east, legend style={draw=none},legend cell align={left},xlabel={$C$}, plotOptions, ymin=0.948, ymax=0.999,xmin=1,xmax=10000,width=0.265\linewidth, height=0.275\linewidth]
		\addplot [color=colorP1] table [x=C,y=a0] {Fig/res_up_light0.999_5.txt};
		\addplot [color=colorP2] table [x=C,y=ae-9] {Fig/res_up_light0.999_5.txt};
		\addplot [color=colorP3] table [x=C,y=ae-8] {Fig/res_up_light0.999_5.txt};
		\addplot [color=colorP4] table [x=C,y=ae-7] {Fig/res_up_light0.999_5.txt};
		\addplot [color=colorP5] table [x=C,y=ae-6] {Fig/res_up_light0.999_5.txt};
		\addplot [color=colorP6] table [x=C,y=ae-5] {Fig/res_up_light0.999_5.txt};
		
		\addlegendentry{$\eta\|\nabla f(x_0)\| = 0$};
		\addlegendentry{$\eta\|\nabla f(x_0)\| = 10^{-6}\tfrac{\mu}{L}$};
		\addlegendentry{$\eta\|\nabla f(x_0)\| = 10^{-5}\tfrac{\mu}{L}$};
		\addlegendentry{$\eta\|\nabla f(x_0)\| = 10^{-4}\tfrac{\mu}{L}$};
		\addlegendentry{$\eta\|\nabla f(x_0)\| = 10^{-3}\tfrac{\mu}{L}$};
		\addlegendentry{$\eta\|\nabla f(x_0)\| = 10^{-2}\tfrac{\mu}{L}$};
		\addplot [color=black,dashed ,domain=1:10000, samples=50] {2*((1-sqrt(0.001))/(1+sqrt(0.001)))^5/(1+((1-sqrt(0.001))/(1+sqrt(0.001)))^10)};
		\addplot [color=black,dotted ,domain=1:10000, samples=50] {0.999^5};
		\addlegendentry{$\rho_*$}
		\addlegendentry{$\rho^k$}
		\end{semilogxaxis}
	\end{tikzpicture}
    \end{tabular}
    \vspace{-0.2cm}
    \caption{Bounds on the convergence rate of Algorithm~\ref{algo:CAA} with $k=5$, $\mu = 10^{-3}$ and $L=1$. Left: $\hat{\rho}(C)$ defined in \eqref{eq:rho-hat-grad}, Middle: right hand side of the bound~\eqref{eq:rho-hat-bound-grad} with $M=k$. Right: right hand side of the bound~\eqref{eq:rho-hat-bound-grad} with $M=1$.}
    \label{fig:example-rate}
\end{figure*}
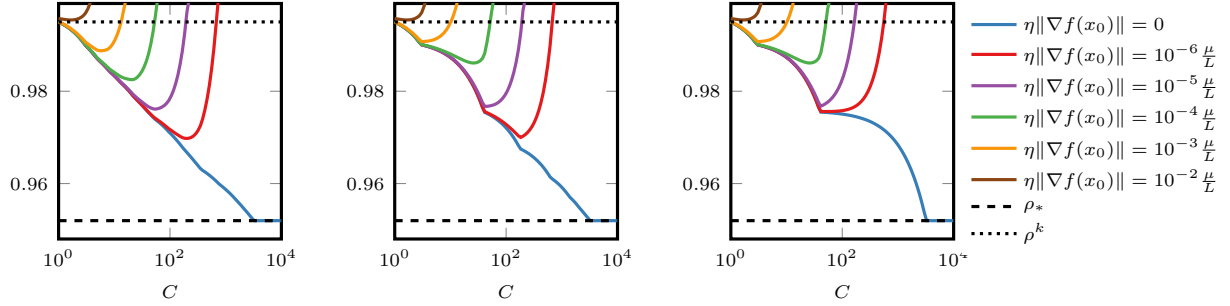
\vspace{-0.3cm}
\begin{remark}
When $\rho \rightarrow 1$,$\alpha_0 \sim \tfrac{1-\rho}{9}$,\\
$\tfrac{2+\rho^k}{2-\rho^k}\rightarrow3$,
$C_1 \rightarrow \tfrac{1}{2}\left(\left(1-\sqrt{2}\right)^k+\left(\sqrt{2}+1\right)^k\right)$,\\ $\alpha_1\sim \min(\tfrac{1}{9},\tfrac{2 (k-1)}{3 \left(\left(1-\sqrt{2}\right)^k+\left(\sqrt{2}+1\right)^k\right)})(1-\rho) $.\\
 $C_* \rightarrow  \tfrac{1}{2}\left(\left(3-2\sqrt{2}\right)^k+\left(2\sqrt{2}+3\right)^k\right)$,\\ $\alpha_2\sim \min(\tfrac{1}{9},\tfrac{2 (2k-1)}{3 \left(\left(3-2\sqrt{2}\right)^k+\left(2\sqrt{2}+3\right)^k\right)})(1-\rho)$.
\end{remark}

\vspace{-0.2cm}
\begin{remark}
Using only the coarse upper bound from Lemma~\ref{lem:cvx-cheb} for $\tilde{\rho}(C)$, we can only guarantee the existence of a $C$ that provides acceleration for Algorithm~\ref{algo:CAA} when $\alpha <  \tfrac{\rho^k-\rho_*}{3kC_*}$. As a comparison,  Proposition~\ref{prop:bound-aa} guarantees acceleration when $\alpha <  \tfrac{\rho^k(1-\rho^k)}{3k(2+\rho^k)}$. When $\rho$ is close to $1$, $\tfrac{\rho^k-\rho_*}{3kC_*}\tfrac{3k(2+\rho^k)}{\rho^k(1-\rho^k)} \sim \tfrac{6 (2 k-1)}{\left(2 \sqrt{2}+3\right)^k+\left(3-2 \sqrt{2}\right)^k}$. 

For instance, when $k=10$ and $\rho$ being close to $1$, our bound guarantees acceleration in Algorithm~\ref{algo:CAA} for $\alpha$'s about $10^6$ times larger than with the coarse bound. 
\end{remark}

\subsection{CAA on Gradient Descent}\label{sec:main-res}
Recently Anderson acceleration has been successfully applied in the optimization field where $F$ is an operator representing an optimization algorithm.

Here we look at the particular case where $F$ is the gradient step operator of a $\mu$-strongly convex function $f : \reals^n \rightarrow \reals$ with $L$-Lipschitz gradient and optimum $x_* \in \reals^n$. In addition, we suppose that $f$ is $\mathcal{C}^2$ with $\eta$-Lipschitz hessian $\nabla^2f$. It is well known (see for instance \cite{Ryu16}) that $F = (I-\tfrac{1}{L}\nabla f)$ is a $\rho=\left(1-\tfrac{\mu}{L}\right)$-Lipschitz operator. We consider Algorithm~\ref{algo:CAA} applied to $F$ at point $x_0\in \reals^n$, and propose to decompose $F$ as 
\[F = \left(I-\tfrac{1}{L}\nabla^2f(x_0)\right) + \tfrac1L\left(\nabla^2f(x_0) - \nabla f\right).\]
We pose $Q = I-\tfrac{1}{L}\nabla^2f(x_0)$ which has its spectrum in $\left[0,\left(1-\tfrac{\mu}{L}\right)\right]$, and $\xi(x) = \tfrac{1}{L}\left(\nabla^2 f(x_0)x - \nabla f(x)\right)$. The following lemma states that on the set where our iterates evolve, $\xi$ is indeed Lipschitz.
\begin{lemma}\label{lem:grad-lip}
Let $k \geq 1$, $C\geq 1$, $x_0\in \reals^n$, $(x_i)_{i\in[1,k]}$ such that $x_{i+1}= x_i - \tfrac{1}{L}\nabla f(x_i)$ for $i\in[0,k-1]$, then $\xi$ is $\tfrac{\eta}{L^2}kC\|\nabla f (x_0)\|$-Lipschitz on $B_C=\{x =\sum_{i=0}^kc_ix_i \in \reals^n|1^tc = 1,\; \|c\|_1 \leq C\}$.
\end{lemma}
\begin{proof} See Appendix~\ref{app:grad-lips}
\end{proof}

\vspace{-0.3cm}
Lemma~\ref{lem:grad-lip} allows applying Proposition~\ref{eq:CAA} to $F$ with $\alpha = \tfrac{\eta}{L^2}kC\|\nabla f(x_0)\|$. There are therefore two ways to have a small $\alpha$ in this context and therefore to guarantee acceleration using Proposition~\ref{prop:bound-aa-tot}, either (i) by having a hessian with a small Lipschitz constant, which means being globally close to a quadratic, or (ii) by being sufficiently close to the optimum (i.e., $\|\nabla f(x_0)\|$ small).

In this setting, $\hat{\rho}(C)$, from Proposition~\ref{prop:bound-aa-tot}, becomes 
\BEQ\label{eq:rho-hat-grad}
\hat{\rho}(C) = \tilde{\rho}(C) + 3\tfrac{\eta}{L^2}\|\nabla f(x_0)\|k^2C^2
\EEQ
and the bound is no more piecewise linear, but piecewise quadratic in $C$, instead.
\BEQ\label{eq:rho-hat-bound-grad}
\begin{aligned}
    \hat{\rho}(C) &\leq \underset{i\in[-1,M]}{\max}\left(\tfrac{C-C_i}{C_{i+1}-C_i}\rho_{i+1} +\tfrac{C_{i+1}-C}{C_{i+1}-C_i}\rho_i,\rho_*\right)\\
    &+3\tfrac{\eta}{L^2}\|\nabla f(x_0)\|k^2C^2
\end{aligned}
\EEQ
As before we study this bound in the case $M=1$ for simplicity. The next proposition provides a range of value of $C$ for which acceleration is guaranteed with Algorithm~\ref{algo:CAA} depending on $\tfrac{\eta}{L^2}\|\nabla f(x_0)\|$
\begin{proposition}\label{prop:bound-aa-grad} Under the assumptions of \S\ref{sec:main-res}.

\begin{enumerate}[label=(\roman*)]
    \item As soon as $\tfrac{\eta}{L^2}\|\nabla f(x_0)\| < \tfrac{\rho^k(1-\rho^k)(2-\rho^k)}{3k^2(2+\rho^k)^2}:= \alpha_3$ there exists  a non empty interval $I$ containing $\tfrac{2+\rho^k}{2-\rho^k}$ such that $\hat{\rho}(C) < \rho^k$ for $C\in I$. 
    \item $\tfrac{\eta}{L^2}\|\nabla f(x_0)\| < \min(\tfrac{\rho^k(1-\rho^k)(2-\rho^k)}{3k^2(2+\rho^k)^2},\tfrac{\rho^k-\rho_1}{3k^2C_1^2}):=\alpha_4$\\
    $\implies \hat{\rho}(C) < \rho^k \text{ for } C\in[\tfrac{2+\rho^k}{2-\rho^k},C_1] $.
\item $\tfrac{\eta}{L^2}\|\nabla f(x_0)\| < \min(\tfrac{\rho^k(1-\rho^k)(2-\rho^k)}{3k^2(2+\rho^k)^2},\tfrac{\rho^k-\rho_*}{3k^2C_*^2}):=\alpha_5$\\
$\implies \hat{\rho}(C) < \rho^k \text{ for } C\in[\tfrac{2+\rho^k}{2-\rho^k},C_*]$.
\end{enumerate}
\end{proposition}
\begin{proof}
Similar to Proposition~\ref{prop:bound-aa}.
\end{proof}
\vspace{-0.2cm}
\begin{remark}
When $\tfrac{\mu}{L} \rightarrow 0$, $\alpha_3 \sim \tfrac{1}{27k}\tfrac{\mu}{L}$.\\
  $\alpha_4\sim \min(\tfrac{1}{27k},\tfrac{4 (k-1)}{3k\left(\left(1-\sqrt{2}\right)^k+\left(\sqrt{2}+1\right)^k\right)^2})\tfrac{\mu}{L} $.\\
 $\alpha_5\sim \min(\tfrac{1}{27k},\tfrac{4 (2 k-1)}{3k\left(\left(2 \sqrt{2}+3\right)^k+\left(3-2 \sqrt{2}\right)^k\right)^2 k})\tfrac{\mu}{L} $.
\end{remark}
\vspace{-0.1cm}
Figure~\ref{fig:example-rate} displays the values of our bounds with fixed $k,\mu$ and $L$ for various values of the perturbation parameter that is $\eta\|\nabla f(x_0)\|$ in this section. We observe that we do not loose much by using the simple upper bound with $M=1$ compared with the exact value of $\tilde{\rho}(C)$ that is to be obtained by solving \eqref{eq:sos}.

Due to the particular form of the perturbation parameter $\alpha$, proportional to $\eta\|\nabla f (x_0)\|$ in the case of the gradient step operator, we see that as soon as $\eta\|\nabla f(x_0)\|$ is small enough to get $\hat{\rho} < 1$, one can apply successively Algorithm~\ref{algo:CAA} while enjoying a perturbation parameter getting smaller and smaller, leading to a faster convergence. Algorithm~\ref{algo:GCAA} is obtained by adding a \emph{guarded step} to the procedure. This step consists in using the extrapolated point $x_e$ only if the gradient at this point is smaller than the gradient at some iterates used in the extrapolation, and allows obtaining nicer convergence properties.
\vspace{-0.1cm}
\begin{algorithm}[!ht]
\caption{Guarded Constrained Anderson Acceleration}
\label{algo:GCAA}
\begin{algorithmic}
\STATE \algorithmicrequire\;$x_0 \in \reals^n$, a strongly function $f$ with $L$-Lipschitz gradient , $C$ a constraint bound and $k$ controlling the number of extrapolation steps, $N$ number of iterations.
\FOR{$i = 0\dots N-1$}
\STATE $x_i^0 = x_i$
\FOR{$j = 0\dots k$}
\STATE $x_{i}^{j+1} = x_i^j - \tfrac{1}{L}\nabla f(x_i^j)$
\ENDFOR
\STATE $R = \begin{bmatrix} x_i^0 -x_i^1&\cdots&x_{i}^k-x_{i}^{k+1}
\end{bmatrix}$ 
\STATE Compute $\tilde{c} = \underset{\substack{\mathbf{1}^Tc = 1,\; \|c\|_1 \leq C}}{\argmin}\;\|Rc\|$ 
\STATE $x_i^e = \sum_{j=0}^{k}\tilde{c}_jx_{i}^j$
\STATE $x_{i+1} = \underset{x \in \{x_i^e,x_i^k\}}{\argmin}\|\nabla f(x)\|$
\ENDFOR
\STATE \algorithmicensure\; $x_{N}$
\end{algorithmic}
\end{algorithm}
\vspace*{-0.3cm}
\begin{proposition}\label{prop:N-iter-online}
Let $f$ be a $\mu$-strongly convex function, with $L$-Lipschitz gradient and $\eta$-Lispchitz hessian. Let $\rho = 1-\tfrac{\mu}{L}$, $k > 2$, $C > 1$, $N \geq 1$ and $x_0 \in \reals^n$, and denote by $(x_i)_{i\in[0,N]}$ the sequence of iterates of Algorithm~\ref{algo:GCAA} on $f$ started at $x_0$ for $N$ iterations with parameter $C$ and $k$. We have
\[ \|\nabla f(x_N)\| \leq \prod_{i=1}^N\hat{\rho}_i(C)\|\nabla f(x_0)\| \]
where 
\begin{align*}
    \hat{\rho}_i & = \min\left(\rho^k,\underset{j\in[-1,1]}{\max}\tfrac{C-C_j}{C_{j+1}-C_j}\rho_{j+1}+\tfrac{C_{j+1}-C}{C_{j+1}-C_j}\rho_j \right.\\
    & + 3\tfrac{\eta}{L^2}\|\nabla f(x_{i-1})\|k^2C^2\Big) 
\end{align*}

In addition,
\[  N \geq \tfrac{\log\left(\tfrac{\eta}{L^2}\tfrac{3k^2(2+\rho^k)^2\|\nabla f(x_0)\|}{\rho^k(1-\rho^k)(2-\rho^k)}\right)}{k\log{\tfrac{1}{\rho}}}
\implies
\prod_{i=1}^N\hat{\rho}_i < \rho^{kN} \]
and $\hat{\rho}_N(C) \underset{N\rightarrow\infty}{\longrightarrow}\underset{j\in[-1,1]}{\max}\tfrac{C-C_j}{C_{j+1}-C_j}\rho_{j+1}+\tfrac{C_{j+1}-C}{C_{j+1}-C_j}\rho_j$.
\end{proposition}
\begin{proof}
This follows by combining results of Proposition~\ref{prop:bound-aa-tot} and Proposition~\ref{prop:bound-aa-grad}.
\end{proof}

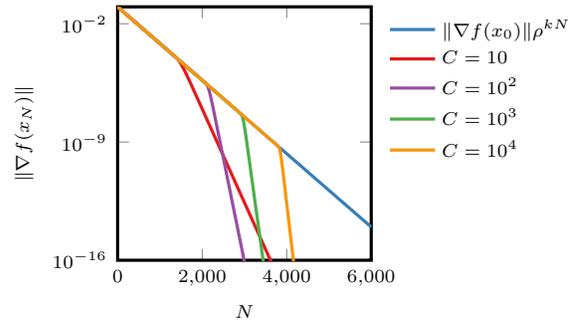
\begin{figure}[!ht]
    \centering
    \begin{tikzpicture}
		\begin{semilogyaxis}[legend pos=outer north east, legend style={draw=none},legend cell align={left},ylabel={$\|\nabla f(x_N)\|$},xlabel={$N$}, plotOptions, ymin=1e-16, ymax=1e-1,xmin=0,xmax=6000,width=0.6\linewidth, height=0.6\linewidth]
		\addplot [color=colorP1] table [x=N,y=r] {Fig/0.999_5_iter_6000.txt};
		\addplot [color=colorP2] table [x=N,y=r10] {Fig/0.999_5_iter_6000.txt};
		\addplot [color=colorP3] table [x=N,y=r100] {Fig/0.999_5_iter_6000.txt};
		\addplot [color=colorP4] table [x=N,y=r1000] {Fig/0.999_5_iter_6000.txt};
		\addplot [color=colorP5] table [x=N,y=r10000] {Fig/0.999_5_iter_6000.txt};
		\addlegendentry{$\|\nabla f(x_0)\|\rho^{kN}$};
		\addlegendentry{$C=10$};
		\addlegendentry{$C=10^2$};
		\addlegendentry{$C=10^3$};
		\addlegendentry{$C=10^4$};
		
		\end{semilogyaxis}
	\end{tikzpicture}
	\vspace{-0.2cm}
    \caption{Illustration of Corollary~\ref{prop:N-iter-online} on with $k=5$, $\mu = 10^{-3}$, $L=1$, $\eta = 10^{-2}$ and $\|\nabla f(x_0)\|=10^{-1}$ for different values of $C$.}
    \label{fig:N-iter-online}
\end{figure}

\vspace{-0.1cm}
The statement of Proposition~\ref{prop:N-iter-online} is illustrated in Figure~\ref{fig:N-iter-online}. In short, as $C$ gets larger, more iterations are required for escaping the guarded regime (as the guarantee from Proposition~\ref{prop:bound-aa-tot} does not provide any improvement over $\rho^k$). But as soon as $\|\nabla f(x_N)\|$ gets small enough, the guarantee of CAA improves that of $F$, escaping the guarded regime, and as $\|\nabla f(x_N)\|$ decreases, the convergence rate gets closer to $\tilde{\rho}(C)$ (which is smaller for large $C$s).

\textbf{Conclusion:} 
We derived upper bounds on the optimal value of the constrained Chebyshev problem, and used them to produce explicit non asymptotic convergence bounds on constrained Anderson acceleration with nonlinear operators. Our convergence bounds are somewhat conservative as they rely on treating the nonlinear part as a perturbation to the linear setting. A remaining open question is whether one can prove better convergence bounds on Anderson acceleration without decoupling linear and nonlinear parts of the operator. However, this would require very different proof techniques.

\clearpage
\onecolumn
\subsubsection*{Acknowledgements}
MB acknowledges support from an AMX fellowship. AT acknowledges support from the European Research Council (grant SEQUOIA 724063). AA is at the d\'epartement d’informatique de l’ENS, \'Ecole normale sup\'erieure, UMR CNRS 8548, PSL Research University, 75005 Paris, France, and INRIA. AA would like to acknowledge support from the {\em ML and Optimisation} joint research initiative with the {\em fonds AXA pour la recherche} and Kamet Ventures, a Google focused award, as well as funding by the French government under management of Agence Nationale de la Recherche as part of the "Investissements d'avenir" program, reference ANR-19-P3IA-0001 (PRAIRIE 3IA Institute).

\bibliographystyle{plainnat}
\bibliography{MainPerso,mybiblio}
\appendix

\section{Proof Proposition~\ref{prop:nonlinear-aa}}\label{app:proof-nonlin-aa}
Let us study the fixed point iterations of $F$, that is
\[ x_{i+1} = G(x_i) + \xi(x_i),\]
which is equivalent to
\[x_{i+1}-x_* = G(x_i-x_*) + \xi(x_i)-\xi(x_*)\]
as $x_* = F(x_*)$. By further developing the previous expression, we arrive to 
\BEQS x_{i+1}-x_\star = G^{i+1}(x_0-x_\star) + \sum_{j=0}^{k}G^{i-j}(\xi(x_j)  -\xi(x_*)).
\EEQS
Another useful quantity is
\BEQS
x_{i+1}- x_i = (G-I)G^i(x_0-x_\star) + (G-I)\sum_{j=0}^{i-1}G^{i-j-1}\left(\xi(x_j)-\xi(x_*)\right)+ \xi(x_i)-\xi(x_*).
\EEQS
Let us use those expressions, along with a triangle inequality, for working out the fixed-point residual
\begin{equation}
    \begin{aligned}
    \lnorm(F-I)(x_e)\rnorm &= \lnorm(G-I)(x_e-x_*) + \xi(x_e)-\xi(x_*)  \rnorm \notag\\
    & = \lnorm\sum_{i=0}^k\tilde{c}_i(G-I)(x_i-x_*) + \xi(x_e)-\xi(x_*) \rnorm\notag\\
    & = \lnorm\sum_{i=0}^k\tilde{c}_iG^i(G-I)(x_0-x_*) +(G-I)\sum_{i=0}^k\tilde{c}_i\sum_{j=0}^{i-1}G^{i-1-j}(\xi(x_j)-\xi(x_*)) + \xi(x_e)-\xi(x_*) \rnorm\notag\\
    & = \lnorm\sum_{i=0}^k\tilde{c}_i(x_{i+1}-x_i) - \sum_{i=0}^k\tilde{c}_i\xi(x_i) - \xi(x_e)\rnorm\notag\\
    & \leq \lnorm\sum_{i=0}^k\tilde{c}_i(x_{i+1}-x_i)\rnorm + \lnorm \xi(x_e) - \sum_{i=0}^k\tilde{c}_i\xi(x_i)\rnorm,\label{ineq:dec_lin_nonlin}
\end{aligned}
\end{equation}
where the first term on the right hand side is exactly the quantity that is minimized in Algorithm~\ref{algo:CAA}. We then bound the two terms separately. Let $c_*$ denotes the coefficients of the polynomial $p_* = \underset{\substack{p \in \reals_k[X]\\ p(1)=1\\\|p\|_1\leq C}}{\argmin}\; \underset{x\in [0,\rho]}{\max} |p(x)|$, we proceed as follows:
\begin{align*}
    &\lnorm\sum_{i=0}^k\tilde{c}_i(x_{i+1}-x_i)\rnorm \leq \lnorm\sum_{i=0}^kc_i^*(x_{i+1}-x_i)\rnorm \text{ since $c_*$ is admissible for the optimization problem \eqref{eq:CAA}}\\
    & = \lnorm\sum_{i=0}^kc_i^*G^i(G-I)(x_0-x_*) +(G-I)\sum_{i=0}^kc_i^*\sum_{j=0}^{i-1}G^{i-1-j}(\xi(x_j)-\xi(x_*)) + \sum_{i=0}^kc_i^*(\xi(x_i)-\xi(x_*)) \rnorm\\
    & = \lnorm\sum_{i=0}^kc_i^*G^i\left[(G-I)(x_0-x_*)+\xi(x_0)-\xi(x_*)\right] +(G-I)\sum_{i=0}^kc_i^*\sum_{j=0}^{i-1}G^{i-1-j}(\xi(x_j)-\xi(x_*)) \right.\\
    &\left.+\sum_{i=0}^kc_i^*\left[\xi(x_i)-\xi(x_*)-G^i(\xi(x_0)-\xi(x_*))\right] \rnorm
\end{align*}
\begin{align*}
    & \leq \lnorm\sum_{i=0}^kc_i^*G^i\left[(F-I)(x_0)\right]\rnorm +\lnorm(G-I)\sum_{i=1}^kc_i^*\sum_{j=0}^{i-1}G^{i-1-j}(\xi(x_j)-\xi(x_*)) \right.\\
    & \left .+ \sum_{i=1}^kc_i^*\left[\xi(x_i)-\xi(x_*)-G^i(\xi(x_0)-\xi(x_*))\right] \rnorm\\
&\leq \lnorm p_*(G)\rnorm\lnorm(F-I)(x_0)\rnorm + \lnorm\sum_{i=1}^kc_i^*\left[\sum_{j=1}^iG^{i-j}(\xi(x_j)-\xi(x_*)) - \sum_{j=0}^{i-1}G^{i-1-j}(\xi(x_j)-\xi(x_*))\right] \rnorm \\
&\leq \lnorm p_*(G)\rnorm\lnorm(F-I)(x_0)\rnorm + \lnorm\sum_{i=1}^kc_i^*\sum_{j=0}^{i-1}G^{i-j-1}\left[\xi(x_{j+1})-\xi(x_j)\right] \rnorm \\
&\leq \lnorm p_*(G)\rnorm\lnorm (F-I)(x_0)\rnorm + \alpha\sum_{i=1}^k|c_i^*|\sum_{j=0}^{i-1}\rho^{i-j-1}\rho^j\lnorm (F-I)(x_0) \rnorm\\
&\leq \left(\lnorm p_*(G)\rnorm + \alpha\sum_{i=1}^k|c_i^*|\rho^{i-1}i\right)\lnorm (F-I)(x_0)\rnorm \\
&\leq \left(\lnorm p_*(G)\rnorm + \alpha k\|c^*\|_1\right)\lnorm (F-I)(x_0)\rnorm \\
&\leq \left(\lnorm p_*(G)\rnorm + \alpha k C\right)\lnorm (F-I)(x_0)\rnorm.
\end{align*}
One can bound $\|p_*(G)\|$ with standard arguments. Since $0 \preccurlyeq G \preccurlyeq \rho I$, there exist an orthogonal matrix $O$ and a diagonal matrix $D$ such that $G = O^tDO$. We get $\|p_*(G)\|=\|O^tp_*(D)O\|\leq\|p_*(D)\|$. One can then notice that $\|p_*(D)\| = \underset{\lambda \in Sp(G)}{\max}|p_*(\lambda)| \leq \underset{x\in[0,\rho]}{\max}|p_*(x)|$, where $Sp(G)$ is the set of eigenvalues of $G$.

Let us bound the second term of the right hand side in~\eqref{ineq:dec_lin_nonlin}
\begin{align*}
    \|\xi(x_e)-\sum_{i=0}^k\tilde{c}_i\xi(x_i)\| & \leq \|\xi(x_e) - \xi(x_{k})\| + \|\xi(x_k) - \sum_{i=0}^k\tilde{c}_i\xi(x_i) \|\\
    & \leq \alpha\left(\|x_e-x_k\| + \sum_{i=0}^k|\tilde{c}_i|\|x_k-x_i\| \right)\\
    & \leq 2\alpha\sum_{i=0}^{k-1}|\tilde{c}_i|\|x_k-x_i\|\\
    & \leq 2\alpha\sum_{i=0}^{k-1}|\tilde{c}_i|\rho^i\|x_{k-i}-x_0\|\\
    & \leq 2\alpha\sum_{i=0}^{k-1}|\tilde{c}_i|\rho^i\sum_{j=0}^{k-1-i}\|x_{j+1}-x_j\|\\
    & \leq 2\alpha\sum_{i=0}^{k-1}|\tilde{c}_i|\rho^i\sum_{j=0}^{k-1-i}\rho^j\|(F-I)(x_0)\|\\
    & \leq 2\alpha\sum_{i=0}^{k-1}|\tilde{c}_i|\rho^i(k-i)\|(F-I)(x_0)\|\\
    & \leq 2\alpha k \|\tilde{c}\|_1\|(F-I)(x_0)\|\\
    & \leq 2\alpha k C\|(F-I)(x_0)\|.
\end{align*}
Combining the two bounds concludes the proof.
\section{Useful Lemmas}\label{app:lemmas}
Unspecified facts on Chebyshev polynomials of the first kind are borrowed from \citep{Maso02}.


\begin{proposition}\label{prop:cheb-pb}
Let $k \in \N$ and $a > 1$, we have
\begin{equation}
    \frac{T_k}{T_k(a)} = \underset{\substack{p\in\reals_k[X]\\p(a)= 1}}{\argmin}\;\underset{x\in[-1,1]}{\max}\;|p(x)|
\end{equation}
where $T_k$ is the first kind Chebyshev polynomials of order $k$. 
\end{proposition}
\begin{proof}
A proof of this result can be found in \citep[Equation 10]{Flan50}. It relies on the fact that $T_k$ reaches its maximum absolute value on $k+1$ points with oscillating sign. 
\end{proof}

Let us show a result that is classically used for analyzing Anderson acceleration~\citep{Golu61,Scie16}.
\begin{proposition}\label{prop:rescaled-cheb}
Let $k \in N$, and $\rho < 1$. 
\[\frac{T_k(\tfrac{2X-\rho}{\rho})}{T_k(\tfrac{2-\rho}{\rho})} = \underset{\substack{p\in\reals_k[X]\\p(1)= 1}}{\argmin}\;\underset{x\in[0,\rho]}{\max}\;|p(x)| \]
where $T_k$ is the first kind Chebychev polynomials of order $k$. And 
\[ \underset{x \in [0,\rho]}{\max}\;\left|\frac{T_k(\tfrac{2X-\rho}{\rho})}{T_k(\tfrac{2-\rho}{\rho})}\right| = \frac{2\beta^k}{1+\beta^{2k}}\]
with $\beta = \tfrac{1-\sqrt{1-\rho}}{1+\sqrt{1-\rho}}$.
\end{proposition}
\begin{proof}
The problem $\underset{\substack{p\in\reals_k[X]\\p(1)= 1}}{\min}\;\underset{x\in[0,\rho]}{\max}\;|p(x)|$ is equivalent to $\underset{\substack{p\in\reals_k[X]\\p(1)= 1}}{\min}\;\underset{y\in[-1,1]}{\max}\;|p(\rho\tfrac{y+1}{2})|$ which is equivalent to $\underset{\substack{q\in\reals_k[X]\\q(\tfrac{2-\rho}{\rho})= 1}}{\min}\;\underset{y\in[-1,1]}{\max}\;|q(y)|$ by denoting $q(y) = p(\rho\tfrac{y+1}{2})$ or equivalently $p(x) = q(\tfrac{2x-\rho}{\rho})$. The last problem is solved using Proposition~\ref{prop:cheb-pb} with $a = \tfrac{2-\rho}{\rho} > 1$. This gives us a solution $q_*(y) = \tfrac{T_k(y)}{T_k(\tfrac{2-\rho}{\rho})}$, and thus a solution to the original problem $p_*(x) = \tfrac{T_k(\tfrac{2x-\rho}{\rho})}{T_k(\tfrac{2-\rho}{\rho})}$.

For the value of the max, we know that $\underset{y \in [-1,1]}{\max}|T_k(y)| = 1$ and then $\underset{x \in [0,\rho]}{\max}\;\left|\tfrac{T_k(\tfrac{2X-\rho}{\rho})}{T_k(\tfrac{2-\rho}{\rho})}\right| = \tfrac{1}{T_k(\tfrac{2-\rho}{\rho})}$. Since $\tfrac{2-\rho}{\rho} > 1$ one can use the formulae for $T_k(x)$ for $|x| \geq 1$ (see e.g \citep[Eq 1.49]{Maso02}): 
\[T_k(x) = \frac12 \bigg( \Big(x-\sqrt{x^2-1} \Big)^k + \Big(x+\sqrt{x^2-1} \Big)^k \bigg) \text{ when } |x| \geq 1.\]
It mechanically follows that
\begin{align*}
    T_k(\tfrac{2-\rho}{\rho}) & = \tfrac{1}{2}  \left( \left(\tfrac{2-\rho}{\rho}-\sqrt{(\tfrac{2-\rho}{\rho})^2-1}\right)^k + \left(\tfrac{2-\rho}{\rho}+\sqrt{(\tfrac{2-\rho}{\rho})^2-1}\right)^k\right)\\
    &= \tfrac{1}{2}  \left( \left(\tfrac{2-\rho-2\sqrt{1-\rho}}{\rho}\right)^k + 
    \left(\tfrac{2-\rho+2\sqrt{1-\rho}}{\rho}\right)^k\right) \\
    & = \tfrac{1}{2}  \frac{\left( \left(1-\sqrt{1-\rho}\right)^{2k} + 
    \left(1+\sqrt{1-\rho}\right)^{2k}\right)}{\rho^k} \\
    & = \tfrac{1}{2}  \frac{\left( \left(1-\sqrt{1-\rho}\right)^{2k} + 
    \left(1+\sqrt{1-\rho}\right)^{2k}\right)}{(1-\sqrt(1-rho)^2)^k} \\
    & = \tfrac{1}{2}  \frac{\left( \left(1-\sqrt{1-\rho}\right)^{2k} + 
    \left(1+\sqrt{1-\rho}\right)^{2k}\right)}{\left(1-\sqrt{1-\rho}\right)^k\left(1+\sqrt{1-\rho}\right)^k} \\ 
    & =  \tfrac{1}{2}\frac{\tfrac{\left(1-\sqrt{1-\rho} \right)^{2k}}{\left(1+\sqrt{1-\rho} \right)^{2k}}+1}{\tfrac{\left(1-\sqrt{1-\rho} \right)^{k}}{\left(1+\sqrt{1-\rho} \right)^{k}}}\\
    & = \tfrac{1+\beta^{2k}}{2\beta^k}
\end{align*}
\end{proof}

The following lemma extends the previous one by looking at polynomials with minimal maximum absolute value on $[-\varepsilon,\rho]$.
\begin{lemma}\label{lem:rescaled-cheb-eps}
Let $k\in \N$, $\rho < 1 $ and $\varepsilon \geq 0$. It holds that
\[ p_\varepsilon = \underset{\substack{p\in\reals_k[X]\\p(1)=1}}{\argmin}\;\underset{x\in[-\varepsilon,\rho]}{\max}\; |p(x)| \]
where $p_\varepsilon(x) = \tfrac{T_k\left(2\tfrac{x+\varepsilon}{\rho + \varepsilon} -1\right)}{\left|T_k\left(2\tfrac{1+\varepsilon}{\rho + \varepsilon} -1\right)\right|} $, and 
\[\underset{x\in[-\varepsilon,\rho]}{\max}|p_\varepsilon(x)| = \tfrac{2\beta_\varepsilon^k}{1+\beta_\varepsilon^{2k}}\]
with $\beta_\varepsilon = \tfrac{1-\sqrt{1-\tfrac{\rho+\varepsilon}{1+\varepsilon}}}{1+\sqrt{1-\tfrac{\rho+\varepsilon}{1+\varepsilon}}}$.
\end{lemma}
\begin{proof} This proof is similar to that of Proposition~\ref{prop:rescaled-cheb}.
Indeed $ \underset{\substack{p\in\reals_k[X]\\p(1)=1}}{\min}\;\underset{x\in[-\varepsilon,\rho]}{\max}\; |p(x)|$ can be reformulated as $\underset{\substack{q\in\reals_k[X]\\q(\tfrac{2(1+\varepsilon)}{\rho+\varepsilon}-1)= 1}}{\min}\;\underset{y\in[-1,1]}{\max}\;|q(y)|$ with $q(y) = p(\tfrac{2(x+\varepsilon)}{\rho+\varepsilon} -1)$. As $\tfrac{2(1+\varepsilon)}{\rho+\varepsilon}-1 \geq 1$ one can apply Proposition~\ref{prop:cheb-pb} and get that $ \underset{\substack{p\in\reals_k[X]\\p(1)=1}}{\min}\;\underset{x\in[-\varepsilon,\rho]}{\max}\; |p(x)|= \tfrac{1}{\left|T_k\left(2\tfrac{1+\varepsilon}{\rho + \varepsilon} -1\right)\right|} = \tfrac{2\beta_\varepsilon^k}{1+\beta_\varepsilon
^{2k}}$ with $\beta_\varepsilon = \tfrac{1-\sqrt{1-\tfrac{\rho+\varepsilon}{1+\varepsilon}}}{1+\sqrt{1-\tfrac{\rho+\varepsilon}{1+\varepsilon}}}$ and $p_\varepsilon(x) = \tfrac{T_k\left(2\tfrac{x+\varepsilon}{\rho + \varepsilon} -1\right)}{\left|T_k\left(2\tfrac{1+\varepsilon}{\rho + \varepsilon} -1\right)\right|} $. 
\end{proof}

In the following, we focus on problems where $\varepsilon$ is close to $0$.

\begin{lemma}\label{lem:rescaled-cheb-oscil} Let $k \in \N^*$ and $\rho < 1$.
For $\varepsilon \in [0,\tilde{\varepsilon}]$ with $\tilde{\varepsilon} = \rho\tfrac{1+\cos(\tfrac{2k-1}{2k}\pi)}{1-\cos(\tfrac{2k-1}{2k}\pi)}$ we have the following properties on $ p_\varepsilon(X)=T_k(\tfrac{2 X}{\rho+\varepsilon} - \tfrac{\rho-\varepsilon}{\rho+\varepsilon})$.
\begin{enumerate}[label=(\roman*)]
\item Let $c \in \reals^{k+1}$ such that $p_\varepsilon(X) = \sum_{i=0}^kc_iX^i$. Then $\mathrm{sign}(c_i) = (-1)^{k-i}$ for $i \in[1,k]$ and  $(-1)^kc_0 \geq 0$.
\item $\left|p_\varepsilon(X)\right|$ is maximal on the $m_i = \tfrac{(\rho+\varepsilon)\cos(\tfrac{i\pi}{k}) + \rho - \varepsilon}{2} \in [-\varepsilon,\rho]$ and $p_\varepsilon(m_i) = (-1)^{i}$.
\end{enumerate}
\end{lemma}
\begin{proof}
The Chebyshev polynomial of first kind $T_k(X)$ is defined such that $T_k(\cos(\theta)) = \cos(k\theta)$ for all $\theta \in \reals$. Thus the roots of $T_k$ are the $(z_i)_{i\in[0 ,k-1]} = (\cos\left(\tfrac{2i+1}{2k}\pi \right))_{i\in[0 ,k-1]} \in [-1,1]$. The roots of $p_\varepsilon(X)$ are the $z^\varepsilon_i$ defined such that $\tfrac{2z^\varepsilon_i}{\rho+\varepsilon} -\tfrac{\rho-\varepsilon}{\rho+\varepsilon} = z_i $. This corresponds to 
\[ z^\varepsilon_i = \tfrac{(\rho+\varepsilon)z_i + \rho - \varepsilon}{2} \in [-\varepsilon,\rho].\]
The smallest root is $z^\varepsilon_{k-1} = \tfrac{(\rho+\varepsilon)\cos(\tfrac{2k-1}{2k}\pi) + \rho - \varepsilon}{2} \geq 0$ for $\varepsilon\in[0,\rho\tfrac{1+\cos(\tfrac{2k-1}{2k}\pi)}{1-\cos(\tfrac{2k-1}{2k}\pi)}]=[0,\tilde{\varepsilon}]$.

This implies that $p_\varepsilon(X)$ \textbf{keeps the same sign on} $[-\varepsilon,0]$. 
In addition $c_0 = 0$ when $\varepsilon=\tilde{\varepsilon}$ and is nonzero otherwise. We also have that $\mathrm{sign}(c_0)=\mathrm{sign}(p_\varepsilon(0))=\mathrm{sign}(p_\varepsilon(-\varepsilon))=\mathrm{sign}(T_k(-1))=(-1)^k$.

For the other coefficients, we first need to show that \textbf{the sign of} $\tfrac{d^lp_\varepsilon}{dx^l}(x)$ \textbf{is also constant on} $[-\varepsilon,0]$. This relies on the Rolle's theorem. Indeed we have $k$ roots in $[0,\rho]$ for $p_\varepsilon$ a polynomial of degree $k$. Thus by Rolle's theorem there exists a root to $p_\varepsilon'(X)$ between each root of $p_\varepsilon$, so the $k-1$ possible roots of $p_\varepsilon'$ are in $]0,\rho[$, and we repeat the argument by applying Rolle's Theorem on $p_\varepsilon'$, etc.

Thus for $i \in [1,k]$, $\mathrm{sign}(c_i)=\mathrm{sign}(\tfrac{d^ip_\varepsilon}{dx^i}(0))=\mathrm{sign}(\tfrac{d^ip_\varepsilon}{dx^i}(-\varepsilon))=\mathrm{sign}(\tfrac{d^iT_k}{dx^i}(-1)) = \mathrm{sign}(\tfrac{d^ip_0}{dx^i}(0))$. Finally we used the formula \citep[Equation 2.12]{Doha2011} for the coefficients of $p_0$ which leads to
\[ \tfrac{d^ip_0}{dx^i}(0) = (-1)^{k-i}k\tfrac{(k+i-1)!}{\Gamma(i+\tfrac{1}{2})(k-i)!\rho^i}\sqrt{\pi}\]
and thus $\mathrm{sign}(c_i) = (-1)^{k-i}$ and we conclude that \textit{(i)} holds.

\textit{(ii)} is a property of the Chebyshev polynomial $T_k$. Indeed since $T_k(\cos(\theta)) = \cos(k\theta)$ then $\underset{x\in[-1,1]}{\max} |T_k(x)| = 1$ and is attained for $x_i= \cos(\tfrac{i\pi}{k})$ with $i\in[0,k]$. In particular $T_k(x_i) = (-1)^i$. Thus $|p_\varepsilon|$ has its maxima on 
\[m_i = \tfrac{(\rho+\varepsilon)\cos(\tfrac{i\pi}{k}) + \rho - \varepsilon}{2} \]
and $p_\varepsilon(m_i) = (-1)^i$.
\end{proof}

\section{Proof of Proposition~\ref{prop:cheb-pb-eps}}\label{app:cheb-pb-eps}

Let show that for $\varepsilon \in [0,\tilde{\varepsilon}]$, $\tilde{\rho}(\|p_\varepsilon\|_1) = \underset{x\in[-\varepsilon,\rho]}{\max}\; |p_\varepsilon(x)| $ with $p_\varepsilon \in \underset{\substack{p\in\reals_k[X]\\p(1)=1}}{\argmin}\;\underset{x\in[-\varepsilon,\rho]}{\max}\; |p(x)|$. 

By Lemma~\ref{lem:rescaled-cheb-eps} we have that $p_\varepsilon$ is a rescaled Chebyshev polynomial. $p_\varepsilon(x) = \tfrac{T_k\left(2\tfrac{x+\varepsilon}{\rho + \varepsilon} -1\right)}{\left|T_k\left(2\tfrac{1+\varepsilon}{\rho + \varepsilon} -1\right)\right|}$. 

\textbf{Goal: } Let show that $p_\varepsilon$ is a local minimum (this will be a global one thanks to convexity) of $p \rightarrow \underset{x\in[0,\rho]}{\max}\;|p(x)|$ on the set $E=\{p \in \reals_k[X],p(1)= 1, \|p\|_1 \leq \|p_\varepsilon\|_1 \}$. 

Let $h = \sum_{i=0}^kh_iX^i \in \reals_k[X] \neq 0$  such that $p_\varepsilon+h \in E$. This implies directly that $h(1) = 0$.

Suppose that 
\BEQ\label{eq:sign-max}
\underset{x\in[0,\rho]}{\max}\;|p_\varepsilon(x)+h(x)| < \underset{x\in[0,\rho]}{\max}\;|p_\varepsilon(x)|.
\EEQ
By definition of $p_\varepsilon$ we have that $\underset{x\in[-\varepsilon,\rho]}{\max}\;|p_\varepsilon(x)+h(x)| \geq \underset{x\in[-\varepsilon,\rho]}{\max}\;|p_\varepsilon(x)|$ and thus by combining it with \eqref{eq:sign-max} we have
\[ \underset{x\in[0,\rho]}{\max}\;|p_\varepsilon(x)+h(x)| < \underset{x\in[0,\rho]}{\max}\;|p_\varepsilon(x)| = \underset{x\in[-\varepsilon,\rho]}{\max}\;|p_\varepsilon(x)| \leq \underset{x\in[-\varepsilon,\rho]}{\max}\;|p_\varepsilon(x)+h(x)|, \]

This implies that  \[\underset{x\in[-\varepsilon,\rho]}{\max}\;|p_\varepsilon(x)+h(x)| = \underset{x\in[-\varepsilon,0]}{\max}\;|p_\varepsilon(x)+h(x)|.\]
Finally this leads to
\BEQ\label{eq:sign-max-1}\underset{x\in[-\varepsilon,0]}{\max}\;|p_\varepsilon(x)+h(x)|  \geq \underset{x\in[-\varepsilon,0]}{\max}\;|p_\varepsilon(x)|.
\EEQ

We write $p_\varepsilon(x) = \sum_{i=0}^k c_iX^i$. By Lemma~\ref{lem:rescaled-cheb-oscil} we know that $\mathrm{sign}(c_i) = (-1)^{k-i}$ and that $|p_\varepsilon|$ is maximal on $[-\varepsilon,\rho]$ at the $m_i =\tfrac{(\rho+\varepsilon)\cos(\tfrac{i\pi}{k}) + \rho - \varepsilon}{2} $ for $i \in [0,k]$ and $\mathrm{sign}(p_\varepsilon(m_i))=(-1)^i$.

In addition, $m_i \in ]0,\rho]$ for $i\in[0,k-1]$. 
Indeed the $m_i$ are in decreasing order and $m_{k-1}$ is strictly larger than the smallest root of $p_\varepsilon$, which means $m_{k-1} > \tfrac{(\rho+\varepsilon)\cos(\tfrac{(2k-1)\pi}{2k}) + \rho - \varepsilon}{2} \geq \tfrac{(\rho+\tilde{\varepsilon})\cos(\tfrac{(2k-1)\pi}{2k}) + \rho - \tilde{\varepsilon}}{2} =0 $ 

Since the $|p_\varepsilon(m_i)|$ are nonzero, we can take $h$ with norm small enough such that $|p_\varepsilon(m_i) + h(m_i)| =|p_\varepsilon(m_i)| + \mathrm{sign}(p_\varepsilon(m_i))h(m_i) = |p_\varepsilon(m_i)| + (-1)^{i}h(m_i)$ for $i \in[0,k-1]$. Equation \eqref{eq:sign-max} imposes $(-1)^ih(m_i) < 0$. Since the $m_i$ are in decreasing order, the mean value theorem tells us that $h$ has a root in each $]m_{i+1},m_i[$ for $i\in[0,k-2]$. In addition, since $p_\varepsilon+h \in E$, $h(1) =0$, meaning that \textbf{$h$ has $k$ roots on $]0,1]$}. Since $h$ is of degree $k$ \textbf{it cannot have any additional root outside $]0,1]$}.


In addition, this means that $h$ \textbf{is nonzero and doesn't change sign on} $]-\infty,0]$. The same conclusion holds for $p_\varepsilon$ on~$]-\infty,0[$.

Taking $\underset{x\in[-\varepsilon/2,0]}{\max}|h(x)| < \underset{x\in[-\varepsilon,-0]}{\max}|p_\varepsilon(x)|-\underset{x\in[-\varepsilon/2,0]}{\max}|p_\varepsilon(x)|$ which is strictly positive since $-\varepsilon$ is the only point of maximum for $p_\varepsilon$ on $[-\varepsilon,0]$, we have that 
\BEQ\max_{x\in[-\varepsilon/2,0]}|p_\varepsilon(x)+h(x)| < \max_{x\in[-\varepsilon,0]}|p_\varepsilon(x)| \leq \max_{x\in[-\varepsilon,0]}|p_\varepsilon(x)+h(x)| \EEQ
by \eqref{eq:sign-max-1}. Thus $|p_\varepsilon+h|$ restricted to $[-\varepsilon,0]$ has its maximum on $[-\varepsilon,-\varepsilon/2]$. 
$|p_\varepsilon| > 0$ and $|h|>0$ on $[-\varepsilon,\tfrac{\varepsilon}{2}]$ and has constant sign on this interval. We denote by $m_* \in [-\varepsilon,-\tfrac{\varepsilon}{2}]$ a point such that $|p_\varepsilon(m_*)+h(x_*)|= \underset{x\in[-\varepsilon,0]}{\max}|p_\varepsilon(x) + h(x)|$. Since $|p_\varepsilon(m_*)| > 0$ one can take $h$ with norm small enough such that $|p_\varepsilon(m_*)+h(m_*)|=|p_\varepsilon(m_*)|+\mathrm{sign}(p_\varepsilon(m_*))h(m_*)= |p_\varepsilon(m_*)|+(-1)^kh(m_*)$. For \eqref{eq:sign-max-1} to hold, $\mathrm{sign}(h(m_*))$ has to be $(-1)^k$. 

Since $h$ has constant sign on $[-\varepsilon,0]$, its sign is $(-1)^k$ and then $\mathrm{sign}(h_0) =\mathrm{sign}(h(0)) =  (-1)^k$.

For $i\in[1,k]$, $c_i\neq 0$ by Lemma~\ref{lem:rescaled-cheb-oscil}, thus we can take $h$ with norm small enough such that $|c_i + h_i| = |c_i| + \mathrm{sign}(c_i)h_i$ for $i\in[1,k]$. Lemma~\ref{lem:rescaled-cheb-oscil} also tells us that $\mathrm{sign}(c_i)=(-1)^{k-i}$ for $i\in[1,k]$ and $(-1)^kc_0\geq 0$. Therefore, one can write \begin{align*}
    \|p_\varepsilon + h\|_1 &= \sum_{i=0}^k|c_i + h_i|\\
    &=  \sum_{i=0}^k(-1)^{k-i}(c_i+h_i) \\
    &= \|p_\varepsilon\|_1 + \sum_{i=0}^k(-1)^{k-i}(h_i)\\
    &= \|p_\varepsilon\|_1 + (-1)^kh(-1)\\
    & \leq \|p_\varepsilon\|_1 \text{ because } p_\varepsilon + h \in E.
\end{align*}  
Finally, we reach 
\[ (-1)^kh(-1) \leq 0 \]
meaning that $\mathrm{sign}(h(-1)) = (-1)^{k+1}$, however we saw that $h$ has sign $(-1)^k$ on $]-\infty,0]$. Since $h$\textbf{ has already all its root  in} $]0,1]$, \textbf{this is a contradiction}.

Thus we cannot find nonzero direction $h$ with arbitrarily small norm such that $p_\varepsilon+h \in E$ and that satisfies~\eqref{eq:sign-max}. $p_\varepsilon$ is thus a local minimum and thus a global one by convexity of the objective and of $E$.

\section{Proof of Proposition~\ref{prop:bound-aa}}\label{app:bound-aa}
Before presenting the proof of the proposition, we show a technical lemma. 
\begin{lemma}\label{lem:order-C}
Let $k\in\N$, $\rho <1$, $C_*$ is defined in \eqref{eq:c*} with explicit value in Remark~\ref{rem:c*} and $C_1$ in \eqref{eq:rho1}. It holds that
\[ \tfrac{2+\rho^k}{2-\rho^k} \leq C_1 \text{ for } k > 1 \]
and 
\[ C_1 \leq C_* \text{ for } k \geq 1 \]
\end{lemma}
\begin{proof}
We start from the expression of $C_1$
\begin{align*}
    C_1 &= \tfrac{\rho_1}{2\rho^k}\left((1-\sqrt{1+\rho^2})^k + (1+\sqrt{1+\rho^2})^k \right)\\
    & = \tfrac{\beta_\rho^k}{\rho^k(1+\beta_\rho^2k)}\left((1-\sqrt{1+\rho^2})^k + (1+\sqrt{1+\rho^2})^k \right),
\end{align*}
where $\beta_\rho = \tfrac{\sqrt{1+\rho}-\sqrt{1-\rho}}{\sqrt{1+\rho}+\sqrt{1-\rho}} = \tfrac{1-1\sqrt{1-\rho^2}}{\rho}$. Thus
\begin{align*}
        C_1 &= \tfrac{(1-\sqrt{1-\rho^2}^k)}{\rho^{2k}+(1-\sqrt{1-\rho^2})^{2k}}\left((1-\sqrt{1+\rho^2})^k + (1+\sqrt{1+\rho^2})^k \right)\\
    &=\tfrac{(1-\sqrt{1-\rho^2}^k)}{(1-\sqrt{1-\rho^2})^k(1+\sqrt{1-\rho^2})^{k}+(1-\sqrt{1-\rho^2})^{2k}}\left((1-\sqrt{1+\rho^2})^k + (1+\sqrt{1+\rho^2})^k \right)\\
    &= \frac{(1-\sqrt{1+\rho^2})^k + (1+\sqrt{1+\rho^2})^k}{(1-\sqrt{1-\rho^2})^k +(1+\sqrt{1-\rho^2})^k}\\
    &= \frac{\sum_{i=0}^{\lfloor k/2 \rfloor} {k\choose 2i}(1+\rho^2)^i}{\sum_{i=0}^{\lfloor k/2 \rfloor} {k\choose 2i}(1-\rho^2)^i}.
\end{align*}
To show that $C_1 \geq  \tfrac{2+\rho^k}{2-\rho^k}$ we need to show
\begin{align*}
    (2-\rho^k)\sum_{i=0}^{\lfloor k/2 \rfloor} {k\choose 2i}(1+\rho^2)^i- (2+\rho^k)\sum_{i=0}^{\lfloor k/2 \rfloor} {k\choose 2i}(1-\rho^2)^i \geq 0,
\end{align*}
and in particular we study 
\[ (2-\rho^k)(1+\rho^2)^i-(2+\rho^k)(1-\rho^2)^i.\]
When $i=0$ this is equal to $-2\rho^k$, when $i=1$ this is equal to $4\rho^2-2\rho^k$. In addition, we can easily see that it is increasing with $i$ and thus it is positive when $i\geq 1$.
We can write for $k \geq 2$
\begin{align*}
    (2-\rho^k)\sum_{i=0}^{\lfloor k/2 \rfloor} {k\choose 2i}(1+\rho^2)^i- (2+\rho^k)\sum_{i=0}^{\lfloor k/2 \rfloor} {k\choose 2i}(1-\rho^2)^i &\geq -2\rho^k + {k\choose 2}(-2\rho^k+4\rho^2)\\
    &\geq 4\rho^2(1-\rho^{k-2})\\
    &\geq 0 \text{ with strict inequality when } k > 2,
\end{align*}
and then $$\boxed{C_1 \geq \tfrac{2+\rho^k}{2-\rho^k}} \text{ with strict inequality when } k > 2.$$

Then we show the second inequality between $C_*$ and $C_1$.
\begin{align*}
    C_* &= \tfrac{\rho_*}{2\rho^k}\left((2+\rho-2\sqrt{1+\rho})^k + (2+\rho+2\sqrt{1+\rho})^k \right)\\
    & = \tfrac{\beta_*^k}{\rho^k(1+\beta_*^2k)}\left((2+\rho-2\sqrt{1+\rho})^k + (2+\rho+2\sqrt{1+\rho})^k \right)
\end{align*}
where $\beta_* = \tfrac{1-\sqrt{1-\rho}}{1+\sqrt{1-\rho}} = \tfrac{2-\rho -2\sqrt{1-\rho}}{\rho}$. Thus $C_*$ can be written as
\begin{align*}
        C_* &= \tfrac{(2-\rho-2\sqrt{1-\rho}^k)}{\rho^{2k}+(2-\rho-2\sqrt{1-\rho})^{2k}}\left((2+\rho-2\sqrt{1+\rho})^k + (2+\rho+2\sqrt{1+\rho})^k \right)\\
    &= \frac{(2+\rho-2\sqrt{1+\rho})^k + (2+\rho+2\sqrt{1+\rho})^k}{(2-\rho-2\sqrt{1-\rho})^k +(2-\rho+2\sqrt{1-\rho})^k}\\
    &= \frac{\sum_{i=0}^{\lfloor k/2 \rfloor} {k\choose 2i}(1+\tfrac{\rho}{2})^{k-2i}(1+\rho)^i}{\sum_{i=0}^{\lfloor k/2 \rfloor} {k\choose 2i}(1-\tfrac{\rho}{2})^{k-2i}(1-\rho)^i}
\end{align*}
when $k\geq 1$
$(1+\tfrac{\rho}{2})^{k-2i}(1+\rho)^i > (1+\rho^2)^i$ and $(1-\tfrac{\rho}{2})^{k-2i}(1-\rho)^i < (1-\rho^2)^i$ for $i\in[0,\lfloor \tfrac{k}{2}\rfloor]$ and thus $$\boxed{C_* > C_1} \text{ when } k \geq 1$$
\end{proof}

\textbf{Proof of Proposition~\ref{prop:bound-aa-tot}:} 
By Lemma~\ref{lem:order-C}, for $k > 2$, $1 < C_0 < C_1 < C_*$.
\begin{enumerate}[label=(\roman*)]
    \item if $\alpha < \tfrac{\rho^k-\rho_0}{3kC_0}=\tfrac{\rho^k-\tfrac{\rho^k}{2-\rho^k}}{3k\tfrac{2+\rho^k}{2-\rho^k}}=\tfrac{\rho^k(1-\rho^k)}{3k(2+\rho^k)}$ then $\hat{\rho}(C_0)= \rho_0 +3\alpha k C_0 < \rho^k$ by Lemma~\ref{lem:tilde-rho-1}.
    \item if $\alpha < \min(\tfrac{\rho^k(1-\rho^k)}{3k(2+\rho^k)},\tfrac{\rho^k-\rho_1}{3kC_1})$ then $\hat{\rho}(C_0) < \rho^k$, and $\hat{\rho}(C_1) \leq \rho_1 +3\alpha k C_1 < \rho^k$. Then by convexity of $\hat{\rho}$, $\hat{\rho}(C) < \rho^k$ for $C \in [C_0,C_1]$. 
    \item if $\alpha < \min(\tfrac{\rho^k(1-\rho^k)}{3k(2+\rho^k)},\tfrac{\rho^k-\rho_*}{3kC_*})$ then $\hat{\rho}(C_0) < \rho_k$ and $\hat{\rho}(C_*) = \rho_* + 3\alpha k C_* < \rho^k$. By convexity of $\hat{\rho}$ we have $\hat{\rho}(C) < \rho^k$ for $C\in[C_0,C_*]$.
\end{enumerate}

\section{Proof of Lemma~\ref{lem:grad-lip}}\label{app:grad-lips}    

We have $D\xi(x) = \tfrac{1}{L}\left(\nabla^2f(x_0). - \nabla^2f(x)\right)$. Under the assumption of $\nabla^2 f$ being $\eta$-Lipschitz, it holds hat $\|D\xi(x)\| \leq \tfrac{\eta}{L}\|x-x_0\|$. Thus for $x\in B_C$
\begin{align*}
    \|D\xi(x)\| &\leq \tfrac{\eta}{L}\|x-x_0\| \\
    &\leq \tfrac{\eta}{L}\|\sum_{i=0}^kc_i(x_i-x_0)\| \\
    &\leq \tfrac{\eta}{L}\sum_{i=0}^k|c_i|\|x_i-x_0\| \\
    &\leq \tfrac{\eta}{L}\sum_{i=1}^k|c_i|\sum_{j=0}^{i-1}\|x_{j+1}-x_j\| \\
    &\leq \tfrac{\eta}{L^2}\sum_{i=1}^k|c_i|\sum_{j=0}^{i-1}\left(1-\tfrac{\mu}{L}\right)^j\|\nabla f(x_0)\| \\
    &\leq \tfrac{\eta}{L^2}k\|c\|_1\|\nabla f(x_0)\|.
\end{align*}
Since $D \xi$ is bounded on the convex set $B_C$, then by the mean value theorem $\xi$ is $\tfrac{\eta}{L^2}kC\|\nabla f(x_0)\|$-Lipschitz on~$B_C$.
\end{document}